\documentclass[a4paper,11pt]{article}

\usepackage{fontawesome5} 
\usepackage{graphicx}

\usepackage{graphicx}
\usepackage{subcaption}
\usepackage{pgfplots}
\pgfplotsset{compat=1.18}

\usepackage{lscape}
\usepackage{longtable}

\usepackage{diagbox}

\usepackage[T1]{fontenc}
\usepackage[latin1]{inputenc}
\usepackage{latexsym} 
\usepackage{amssymb}
\usepackage{fancyvrb}
\usepackage{algorithmic}

\usepackage[T1]{fontenc}
\usepackage[latin1]{inputenc}
\usepackage[english]{babel}
\usepackage{latexsym} 
\usepackage{amssymb}
\usepackage{amsmath}
\usepackage{array}
\usepackage{bbding} 
\usepackage{mathpazo} 

\usepackage{amsthm}
\usepackage{xcolor}

\usepackage{apacite}

\usepackage{setspace}
\usepackage[font=small]{caption}
\usepackage{rotating}

\usepackage{subcaption}
\usepackage{tikz}
\usetikzlibrary{backgrounds}

\usetikzlibrary{shapes.geometric}
\usetikzlibrary{arrows,automata}

\usepackage{multirow}

\usepackage{color}
\usepackage{colortbl}

\usepackage{pgf}
\usepackage{mathrsfs}
\usetikzlibrary{arrows}

\usepackage{algorithmic}
\usepackage{longtable}
\usepackage{arydshln}
\usepackage{hyperref}
\usepackage{pdflscape}

\usepackage{listings}

\allowdisplaybreaks
\usepackage[margin=0.8 in]{geometry}
\usepackage{natbib}
\PassOptionsToPackage{numbers, compress}{natbib}

\theoremstyle{plain}

\newtheorem{theorem}{Theorem}[section]
\newtheorem{proposition}[theorem]{Proposition}

\newtheorem{definition}[theorem]{Definition}
\newtheorem{example}[theorem]{Example}

\theoremstyle{definition}

\theoremstyle{remark}

\usepackage{floatrow}
\newfloatcommand{capbtabbox}{table}[][\FBwidth]

\usepackage{blindtext}

\newcommand{\addQEDstyle}[2]{\AtBeginEnvironment{#1}{\pushQED{\qed}\renewcommand{\qedsymbol}{#2}}\AtEndEnvironment{#1}{\popQED}}
\addQEDstyle{example}{$\triangle$}

\usepackage{siunitx}

\newcolumntype{R}[2]{%
    >{\adjustbox{angle=#1,lap=\width-(#2)}\bgroup}%
    l%
    <{\egroup}%
}
\allowdisplaybreaks

\usepackage{pgfplots}
\usepackage{pgfplotstable}
\usepackage{subcaption}
\usepackage{float}
\pgfplotsset{compat=1.18}

\pgfplotstableread[col sep=space]{
i prop cea cel ep
1 6131.527 6128.253 6132.778 6130.500
2 6221.583 6218.286 6222.811 6220.566
3 6233.590 6230.290 6234.816 6232.575
4 6239.594 6236.293 6240.818 6238.580
5 6503.758 6500.390 6504.915 6502.774
6 6665.858 6662.450 6666.975 6664.893
7 6773.924 6770.489 6775.015 6772.972
8 7158.158 7154.631 7159.156 7157.254
9 7488.356 7484.752 7489.278 7487.496
10 7512.371 7508.761 7513.286 7511.514
11 7620.435 7616.801 7621.326 7619.593
12 7794.539 7790.865 7795.390 7793.721
13 7956.635 7952.925 7957.450 7955.840
14 7980.649 7976.934 7981.459 7979.858
15 8058.695 8054.962 8059.488 8057.915
16 8250.809 8247.033 8251.558 8250.056
17 8268.819 8265.040 8269.565 8268.069
18 8292.833 8289.049 8293.574 8292.087
19 8611.019 8607.166 8611.691 8610.321
20 8791.124 8787.232 8791.757 8790.453
21 8827.145 8823.245 8827.771 8826.479
22 8911.194 8907.276 8911.802 8910.541
23 8995.242 8991.307 8995.833 8994.603
24 9157.336 9153.367 9157.892 9156.722
25 9193.356 9189.380 9193.906 9192.748
26 9223.373 9219.391 9223.917 9222.770
27 9325.432 9321.429 9325.954 9324.845
28 9337.439 9333.433 9337.959 9336.854
29 9493.528 9489.491 9494.016 9492.969
30 9583.580 9579.524 9584.049 9583.035
31 9721.658 9717.575 9722.100 9721.136
32 9853.733 9849.624 9854.149 9853.233
33 10075.859 10071.706 10076.231 10075.396
34 10376.028 10371.816 10376.341 10375.617
35 10400.042 10395.825 10400.350 10399.634
36 10808.270 10803.976 10808.501 10807.934
37 10862.300 10857.996 10862.521 10861.974
38 11060.410 11056.069 11060.594 11060.119
39 11114.440 11110.089 11114.614 11114.159
40 11246.514 11242.137 11246.662 11246.256
41 11270.527 11266.146 11270.671 11270.273
42 11336.563 11332.170 11336.696 11336.322
43 11582.699 11578.261 11582.786 11582.502
44 11588.703 11584.263 11588.789 11588.507
45 11852.848 11848.361 11852.886 11852.701
46 11900.874 11896.379 11900.904 11900.736
47 11978.917 11974.407 11978.933 11978.793
48 12104.986 12100.454 12104.979 12104.886
49 12321.104 12316.534 12321.059 12321.045
50 12417.156 12412.569 12417.094 12417.115
51 12771.348 12767.699 12771.225 12771.375
52 12777.351 12773.720 12777.227 12777.380
53 12807.368 12803.829 12807.238 12807.402
54 12873.403 12870.066 12873.263 12873.450
55 12993.468 12990.497 12993.307 12993.538
56 13203.581 13201.247 13203.384 13203.693
57 13269.617 13267.481 13269.409 13269.741
58 13335.652 13333.715 13335.433 13335.790
59 13809.906 13809.376 13809.608 13810.138
60 13827.916 13827.439 13827.615 13828.151
61 13971.993 13971.938 13971.668 13972.257
62 13990.002 13990.000 13989.675 13990.270
63 14104.063 14104.393 14103.717 14104.354
64 14218.124 14218.784 14217.759 14218.437
65 14374.207 14375.317 14373.817 14374.552
66 14488.267 14489.704 14487.859 14488.636
67 14560.305 14561.948 14559.885 14560.689
68 14572.312 14573.989 14571.890 14572.697
69 14680.369 14682.354 14679.930 14680.777
70 15100.590 15103.759 15100.085 15101.085
71 15112.597 15115.799 15112.089 15113.094
72 15184.635 15188.038 15184.116 15185.147
73 15202.644 15206.098 15202.123 15203.160
74 15298.694 15302.415 15298.158 15299.231
75 15418.757 15422.811 15418.202 15419.319
76 15694.902 15699.714 15694.304 15695.522
77 15754.933 15759.910 15754.327 15755.566
78 15838.977 15844.183 15838.358 15839.628
79 15844.980 15850.202 15844.360 15845.632
80 16181.156 16187.287 16180.484 16181.879
81 16211.171 16217.383 16210.495 16211.901
82 16253.193 16259.517 16252.511 16253.932
83 16457.299 16464.169 16456.586 16458.082
84 16499.321 16506.303 16498.602 16500.112
85 16523.333 16530.379 16522.611 16524.130
86 16601.374 16608.626 16600.639 16602.187
87 16751.452 16759.101 16750.695 16752.297
88 16799.476 16807.252 16798.713 16800.333
89 16895.526 16903.554 16894.748 16896.403
90 16961.560 16969.761 16960.773 16962.452
91 16997.579 17005.874 16996.786 16998.478
92 17165.665 17174.399 17164.848 17166.602
93 17177.672 17186.436 17176.852 17178.610
94 17195.681 17204.492 17194.859 17196.624
95 17375.774 17385.052 17374.926 17376.756
96 17393.783 17403.107 17392.932 17394.769
97 17693.937 17704.032 17693.043 17694.990
98 17789.986 17800.327 17789.079 17791.060
99 17850.017 17860.510 17849.101 17851.104
100 17916.051 17926.712 17915.125 17917.153
}\allocations

\pgfplotstablecreatecol[
    create col/expr={\thisrow{prop}-\thisrow{cea}}
]{diffPC}\allocations

\pgfplotstablecreatecol[
    create col/expr={\thisrow{prop}-\thisrow{cel}}
]{diffPT}\allocations

\pgfplotstablecreatecol[
    create col/expr={\thisrow{cea}-\thisrow{cel}}
]{diffCT}\allocations

\pgfplotstablecreatecol[
    create col/expr={\thisrow{ep}-\thisrow{prop}}
]{diffEPprop}\allocations

\pgfplotstablecreatecol[
    create col/expr={\thisrow{ep}-\thisrow{cea}}
]{diffEPcea}\allocations

\pgfplotstablecreatecol[
    create col/expr={\thisrow{ep}-\thisrow{cel}}
]{diffEPcel}\allocations

\begin{document}

\title{\textbf{Group Buying through Bankruptcy Rules: Cost Sharing and Endogenous Pricing}}

\author{M. G. Fiestras-Janeiro$^1$, A. Meca$^2$, A. Saavedra-Nieves$^{\ast,3}$}

\date{\footnotesize \emph{
		$^1$ CITMAga. University of Vigo. Departament of Statistics and Operations Research. \texttt{fiestras@uvigo.gal}  \\
        $^2$ i-MInD (Instituto Universitario de Investigaci\'on en Matem\'aticas, 
Ciencia de Datos e Inteligencia Artificial), Universidad  Miguel Hern\'andez de Elche, Elche, Spain. \texttt{ana.meca@umh.es}\\ 
		$^3$  \emph{Corresponding author.} CITMAga. University of Santiago de Compostela. Department of Statistics, Mathematical Analysis, and Optimization.  \texttt{alejandro.saavedra.nieves@usc.es} \\}}

\maketitle

\begin{abstract}\color{black}
Group buying allows independent buyers to pool their demand and benefit from quantity discounts that would not be available to them individually. While collective purchasing reduces total procurement costs, it raises a fundamental allocation question: how should the resulting savings be distributed among heterogeneous buyers? We address this question by integrating cost-allocation and bankruptcy approaches into a unified group-buying framework. Starting from a continuous quantity discount function, we examine cost-gap-based allocation mechanisms and develop an endogenous pricing representation in which each buyer's effective discount may depend on the demand profile of the entire group. We then establish a connection between group buying and bankruptcy theory by interpreting the aggregate savings from cooperation as the estate and each buyer's potential savings as its claim. This formulation allows us to derive cost-allocation mechanisms based on the proportional, constrained equal awards, constrained equal losses, and Talmud rules. For each rule, we obtain the corresponding buyer-specific price function and characterize the associated variable discount factor. We further examine the relationships between bankruptcy-based and cost-gap-based allocation mechanisms and use a numerical example with simulated data to compare their implications. The proposed framework links alternative fairness principles for allocating cooperative savings to economically interpretable pricing mechanisms, providing a unified basis for analyzing how group purchasing gains are distributed among heterogeneous buyers.
\end{abstract}
\vspace{0.25cm}

\noindent {\bf Keywords:} group buying, cost allocation, bankruptcy rules, endogenous pricing, quantity discounts

\section{Introduction}\label{sec1}

Group purchasing has become an increasingly relevant procurement strategy for firms and organizations seeking to exploit economies of scale. By pooling their demands, otherwise independent buyers may reach purchase volumes that allow them to obtain prices that would not be available individually. In practice, group purchasing is used in a variety of settings, including retail supply chains \citep{zhou2014potentially}, small and medium-sized enterprises (SMEs) purchasing groups \citep{ghaderi2012empirical}, healthcare procurement \citep{nollet2003development,toulemon2018effect}, and online markets \citep{zhang2016optimal}. These applications illustrate the relevance of demand aggregation across commercial supply chains, institutional procurement, and digital platforms. Their economic rationale is straightforward: when unit prices decrease with purchased quantity, demand aggregation can generate substantial savings for the participating buyers.


The literature on cooperative purchasing and group buying has examined this phenomenon from several complementary perspectives. Early contributions emphasized the organizational design of cooperative purchasing arrangements and the mechanisms through which buyers can coordinate their purchases to exploit quantity discounts. \cite{schotanus2007developing} analyze different forms of cooperative purchasing and show that the organizational structure of the purchasing group affects its advantages, disadvantages, and critical success factors. From a mechanism-design perspective, \cite{li2004mechanism} study coalition formation and cost sharing in group-buying markets, focusing on properties such as efficiency, budget balance, individual rationality, incentive compatibility, and coalition stability. \cite{chen2011group} further show that the benefits of group buying need not be distributed uniformly among participants, since buyer heterogeneity can substantially affect the attractiveness of cooperation. More recent research has extended these ideas by incorporating strategic interactions, endogenous pricing mechanisms, and alternative procurement arrangements. \cite{guan2024group}, for instance, compare group-buying pricing with flat pricing and individual quantity discounts in a business-to-business supply chain with multiple retailers, highlighting the role of buyer heterogeneity. \cite{fu2024group} study group buying among competing firms under alternative power structures, while \cite{du2025does} analyze group buying and agency procurement as alternative forms of joint procurement under quantity discounts. 
More recently, \cite{zhang2026fresh} analyze a fresh-products community group-buying delivery problem with heterogeneous customers, explicitly incorporating differences in customer behavior and pickup-time preferences. Their results show that customer heterogeneity can significantly affect the operational performance of group-buying systems. Although their focus is on delivery optimization rather than cost allocation, this evidence further supports the importance of accounting for heterogeneity when designing collective purchasing arrangements. Taken together, these studies show that group buying remains an active research area in which buyer heterogeneity, procurement design, and the distribution of cooperative gains are important considerations.

A closely related stream of literature concerns the purchasing and pricing mechanisms underlying procurement decisions. \cite{schotanus2009unraveling} propose a continuous quantity discount function that provides a parsimonious representation of the relationship between purchased quantity and unit price and can be incorporated into cooperative purchasing models. Importantly for the present study, they identify the calculation and allocation of savings generated by pooled demand as a relevant problem for purchasing groups. Related research has examined procurement design and purchasing decisions in broader operational settings. \cite{arantes2022application} develop a purchasing portfolio approach for defining medicine procurement strategies, emphasizing the strategic role of purchasing decisions, while \cite{meng2023emergency} design a bidirectional option-contract mechanism for emergency-supply procurement and study how contractual terms can coordinate purchasing decisions under uncertainty. More recently, \cite{rahiminia2026design} incorporate quantity discounts into sustainable supply-chain design under supply-side disruption. Although these studies address different procurement settings and do not directly focus on the allocation of cooperative savings in group buying, they illustrate the broader relevance of purchasing and pricing mechanisms in determining procurement outcomes.


The existence of aggregate savings, however, does not determine how those savings should be distributed. Buyers may differ substantially in their individual demands and, consequently, in their potential contribution to the price reductions achieved by the group. Buyers with large demands may emphasize their role in reaching favorable purchase volumes, whereas smaller buyers also contribute to aggregate demand and may require sufficient incentives to participate. The allocation problem therefore involves a choice among different distributional principles, including proportionality, equal treatment, and the distribution of gains according to buyers' claims. These principles may lead to different individual payments even when the total purchasing cost remains unchanged.


One approach to this problem relies on cost-allocation methods. In particular, the cost-gap methodology (cf. \citealp{driessen1985cost}) decomposes total cost into separable costs, associated with individual marginal contributions, and a non-separable component generated by joint participation. In the present group-buying setting, this decomposition makes it possible to distinguish the portion of the purchasing cost attributable to each buyer from the interaction effects created by demand aggregation. It also provides a basis for distributing the non-separable component according to alternative principles. In our model, these principles give rise to the Egalitarian Non-Separable Cost (ENSC) and Surplus Cost Ratio Based (SCRB) rules, which can be compared with the Equal Price (EP) allocation (cf. \citealp{schotanus2008unfair}).


Although previous studies have addressed coalition formation, cost sharing, stability, and incentive compatibility in group-buying environments, comparatively less attention has been paid to the systematic allocation of cooperative savings among heterogeneous buyers using alternative normative principles. In particular, the potential connection between group-buying cost allocation and bankruptcy theory has not been systematically explored. This connection is especially natural because collective purchasing generates a finite amount of savings to be divided among buyers whose potential claims on those savings may exceed the amount available. This structure closely resembles a bankruptcy problem.


In the classical formulation, a bankruptcy problem consists of an estate to be divided among claimants whose aggregate claims exceed the estate. Bankruptcy theory therefore does not merely provide an alternative collection of allocation formulas; it offers a conceptually different way of interpreting the allocation of cooperative gains. The bankruptcy literature has developed a rich family of division rules based on distinct normative principles, including proportionality, equal awards, equal losses, and combinations of these ideas. O'Neill's seminal analysis already emphasized the role of claims and alternative principles for dividing an insufficient estate (cf. \citealp{o1982problem}). The bankruptcy literature has developed a rich family of division rules based on distinct normative principles, including proportionality, equal awards, equal losses, and combinations of these ideas (cf. \citealp{thomson2003axiomatic}; \citealp{thomson2015axiomatic}).

The main idea of this paper is that a bankruptcy problem need not be imposed exogenously on a group-buying situation. Rather, it is generated endogenously by the cooperative purchasing technology itself. For every buyer, we define a claim from the difference between the cost of purchasing independently and the buyer's separable cost within the group. The corresponding estate is the aggregate saving obtained when the buyers purchase jointly rather than individually. Under the quantity discount structure considered here, these quantities satisfy the requirements of a well-defined bankruptcy problem. The cooperative purchasing problem can therefore be transformed into a claims problem in which classical bankruptcy rules determine how the gains from cooperation are shared.

This bankruptcy interpretation also provides a useful connection between cost allocation and pricing. Standard quantity discount functions typically express the unit price as a function of the quantity purchased by a buyer. In a collective purchasing environment, however, the effective price assigned to one buyer is naturally affected by the demands of the other participants. The composition of the group matters because it determines both the total discount generated by cooperation and its distribution among members. We therefore introduce a variable-discount formulation in which the unit price of buyer \(i\) depends on both its own demand \(q_i\) and the demand vector of the remaining buyers \(q_{-i}\). This formulation allows each allocation rule to be represented by an induced price function, thereby translating alternative principles for sharing cooperative gains into alternative pricing schemes.

Against this background, the paper makes several contributions. First, we consider a group-buying model based on the continuous quantity discount function in \cite{schotanus2009unraveling} and assess well-known allocations in literature derived from cost-gap principles (cf. \citealp{driessen1985cost}). 
Second, we introduce a general endogenous-pricing representation that allows the discount parameter faced by an individual buyer to depend on the entire demand profile. Third, we establish the endogenous bankruptcy structure associated with group buying and use it to construct cost allocations based on four classical bankruptcy rules: the proportional rule, the constrained equal awards rule, the constrained equal losses rule, and the Talmud rule. For each of them, we derive the corresponding bankruptcy-driven group-buying price function and show how cooperation modifies the original discount factor. Fourth, we compare the bankruptcy-driven allocations with cost-gap-based mechanisms, identifying cases in which the two approaches coincide as well as situations in which they distribute cooperative gains differently. Finally, a numerical analysis with heterogeneous buyers illustrates the behavior of the proposed mechanisms and the distributional consequences of alternative allocation principles.

An important feature of the proposed framework is that it separates the generation of cooperative savings from the principle used to distribute them. The quantity discount function determines the economic gains generated by pooling demand, whereas the allocation rule determines how these gains are translated into buyer-specific payments and prices. Consequently, different notions of fairness can be incorporated without changing the underlying purchasing technology. This separation provides both analytical flexibility and managerial interpretability: the choice of a cost-sharing principle can be directly understood in terms of the effective discount received by each participant.


The remainder of the paper is organized as follows. Section \ref{sec2} introduces the group-buying model, develops the cost-gap decomposition, and presents the variable-discount pricing framework. Section \ref{sec3} constructs the bankruptcy problem associated with group buying, derives the allocations and induced price functions under the four bankruptcy rules, and examines their relationships with cost-gap-based mechanisms. Section \ref{sec4} presents a numerical example and a simulation study to investigate the effects of demand heterogeneity, group size, and the parameter regimes considered on allocation outcomes and variable discount factors. Section \ref{sec5} concludes and discusses directions for future research. Appendices \ref{app0}-\ref{app3} includes some extra results.

\section{The Group Buying Model}\label{sec2}

In this section, we describe the multi-agent group buying situation considered as a basis of this paper.

We consider a finite set of agents $N=\{1,\ldots,n\}$, who want to buy a certain service or good from a supplier. Each agent $i\in N$ needs a certain amount $q_i>0$ of units of good at a cost $P(q_i)q_i$. In this context, we consider a function $P:\mathbb{R}_{+}\rightarrow \mathbb{R}$ that represents the variable cost function per unit in such a way that $P(q_i)>0$ denotes the price that agent $i$ pays for quantity $q_i$. We assume that $P$ is a non-increasing function and that $P(q_i)q_i$ (the price function) is a concave function. A common example of the function $P$ is the standard price function, originally introduced in \cite{schotanus2009unraveling} and defined for all $q\in \mathbb{R}^{+}$, by
	\begin{equation}\label{eq:psch}
		P(q)=a+\frac{b}{q^{\eta}} \mbox{ for all }q\in \mathbb{R}^{+}
	\end{equation}
	where $a\geq 0$,  $b$ is a discount factor (positive or negative) and $\eta$ is the parameter that represents the steepness of a quantity discount function. To ensure that $P$ is well defined when $b<0$, its domain is assumed to be a compact real interval, while $b\eta>0$ guarantees that $P$ is non-increasing and $\eta\in[-1,1]$ ensures concavity. We restrict our study to the case $b>0$ and $\eta\in(0,1)$.

 We refer to our model as a group-buying multi-agent situation $(N,q,P)$ (henceforth, GB-model), where $N$ is the set of agents involving in the group-buying situation, $q=(q_i)_{i\in N}$ denotes the vector of quantities of goods that each $i\in N$ needs, and $P$ is a cost function under the above-mentioned assumptions. The set of group-buying multi-agent situation is denoted by $GB^N$.

\subsection{A cost gap approach to group-buying allocation}\label{costgap}

A relevant issue in group-buying environments is the allocation of the total cost generated by collective purchasing among the participating agents. For any $(N,q,P)\in GB^N$, the total cost of the group is given by
\[
C^{GB}(N,q,P)=P(q_N)q_N,
\]
where $q_N=\sum_{i\in N}q_i$ denotes the aggregate demand of the group and $P(q_N)$ is the  price obtained under collective
purchasing by the group. In the case of a single agent $i$, the individual cost is given by $C^{GB}(\{i\},q_i,P)=P(q_i)q_i$.

An allocation rule in group-buying is a mapping $\phi$ that assigns to each group-buying situation $(N,q,P)$ an allocation vector $\phi(N,q,P)\in\mathbb{R}^{|N|}$ such that efficiency holds, that is,
\[
\sum_{i\in N}\phi_i(N,q,P)=C^{GB}(N,q,P).
\]

A natural benchmark in this context is the Equal Price rule, which assigns to each agent a proportional share of the total cost according to her individual demand. Formally, it is given by
\[
EP_i(N,q,P)=P(q_N)q_i.
\]
Under the Schotanus-type price function in Equation (\ref{eq:psch}), as $P(q_N)=a+\frac{b}{q_N^{\eta}}$, 
it holds that
\[
EP_i(N,q,P)=aq_i+b\frac{q_i}{q_N^{\eta}},
\]
which highlights that each agent pays a linear component plus a decreasing term that reflects economies of scale generated by the group.

Following the cost gap methodology in cooperative cost allocation (cf. \citealp{driessen1985cost}), a key idea is that total costs can be decomposed into a separable component, which captures marginal contributions of agents, and a non-separable component, which captures interaction effects arising from joint participation.

In this spirit, the \emph{separable cost of agent} $i$ as
\begin{align}\label{SCrule}
    SC_i(N,q,P)&
=
C^{GB}(N,q,P)-C^{GB}(N\setminus\{i\},q_{-i},P)\notag\\
&=q_iP(q_N)
-
\big(P(q_{N\setminus \{i\}})-P(q_N)\big)q_{N\setminus \{i\}}.
\end{align}The first term, $q_iP(q_N)$, represents the cost that agent $i$ would face under the prevailing group price, while the second term captures the change in total cost induced by the removal of agent $i$, thus reflecting her marginal contribution to the economies of scale. Such expression can be explicitly written, for the price function in Equation (\ref{eq:psch}),  as
\[
SC_i(N,q,P)
=
aq_i+b(q_N^{1-\eta}-q_{N\setminus \{i\}}^{1-\eta}).
\]
\noindent Under the Schotanus price specification, the separable cost of agent $i$ can be expressed as the Equal Price allocation minus the reduction in the aggregate cost of the remaining agents in the group. We may notice that $C^{GB}(\{i\},q_{-i},P)>SC_i(N,q,P)$, for all $i\in N$, since the function $f(x)=x^{1-\eta}$ is strictly concave (see details in Appendix \ref{appNSC}). Moreover, $\sum_{i\in N}C^{GB}(\{i\},q_{-i},P)=C^{GB}(N,q,P).$

Once separable costs are assigned, the remaining part of the total cost defines the non-separable cost component, given by
\begin{align}NSC(N,q,P)
=&
C^{GB}(N,q,P)-\sum_{i\in N}SC_i(N,q,P)\\
   =&
\sum_{i\in N}
\big(P(q_{N\setminus \{i\}})-P(q_N)\big)q_{N\setminus \{i\}}, 
\end{align}
that is strictly non-negative when the price function $P$ is a non-increasing function.

 When using the price function in Equation (\ref{eq:psch}), this expression becomes
\begin{equation}\label{NSCformula}
    \begin{aligned}
 NSC(N,q,P)
 =& b\bigg(q_N^{1-\eta}-\sum_{i\in N}(q_N^{1-\eta}-q_{N\setminus \{i\}}^{1-\eta}  ) \bigg),\\
 \end{aligned}
\end{equation}being $NSC(N,q,P)>0$ (see details in Appendix \ref{appNSC}). {\color{black}Hence, $\sum_{i\in N}SC_i(N,q,P)<C^{GB}(N,q,P)$.}

This term represents the cost differential between each reduced coalition and the grand coalition, weighted by the corresponding aggregate demand. Therefore, it captures the portion of total costs that cannot be directly associated with separable individual contributions.

To allocate this interaction component, we present an egalitarian principle and consider the Egalitarian Non-Separable Cost (ENSC) rule  for group buying, which assigns to each agent an equal share of the non-separable cost. Formally,
\begin{equation}\label{ENSCrule0}
ENSC_i(N,q,P)
=
SC_i(N,q,P)
+
\frac{1}{n}NSC(N,q,P).
\end{equation}

This rule combines individual marginal responsibility (through $SC_i$) with an equal redistribution of the collective interaction effects (through $NSC$). In this sense, it embodies a two-layer allocation principle: first, agents are charged according to their marginal impact on total cost, and second, the residual cooperation surplus or deficit is shared equally among all participants. The expression of the ENSC rule, when considering the price function given in (\ref{eq:psch}) and  since $\sum_{j\in N}q_{N\setminus \{j\}}=(n-1)q_N$, is given by
\begin{equation}\label{ENSCrule}
    \begin{aligned}
ENSC_i(N,q,P)&=aq_i+b(q_N^{1-\eta}-q_{N\setminus \{i\}}^{1-\eta})+\frac{b}{n}\bigg(q_N^{1-\eta}-\sum_{j\in N}(q_N^{1-\eta}-q_{N\setminus \{j\}}^{1-\eta}  ) \bigg). \\ 
&= aq_i+\frac{b}{n}\bigg(q_N^{1-\eta}+\sum_{j\in N}(q_{N\setminus \{j\}}^{1-\eta}-q_{N\setminus \{i\}}^{1-\eta})\bigg)
\end{aligned}
\end{equation}for every $i\in N$. From an economic perspective, the ENSC rule provides a decomposition of group-buying costs into individual and collective components. The separable part reflects the marginal cost contribution of each agent under the group price, while the non-separable part captures the interaction generated by joint purchasing decisions. The egalitarian allocation of this latter component ensures that all agents equally share the effects of cooperation that cannot be attributed to individual behavior, thereby introducing a natural notion of fairness in the distribution of group-generated cost externalities.

A natural alternative to the egalitarian treatment of the non-separable cost component is to distribute it proportionally to the deviation between individual costs and separable costs. This leads to the so-called Surplus Cost Ratio Based (SCRB) method. Formally, for any $(N,q,P)\in GB^N$, the SCRB rule for group buying is defined, for every $i\in N$, by
\begin{equation}\label{SCRBrule}
    SCRB_i(N,q,P)
=
SC_i(N,q,P)
+
\frac{P(q_i)q_i - SC_i(N,q,P)}
{\sum_{j\in N}\big(P(q_j)q_j - SC_j(N,q,P)\big)}
\cdot NSC(N,q,P).
\end{equation}

\noindent This rule assigns to each agent his separable cost component and distributes the non-separable cost proportionally to the gap between his individual cost and his separable cost. In particular, the weight
\[
P(q_i)q_i - SC_i(N,q,P)
\]
measures how much agent $i$ contributes to the difference between his isolated purchasing and his marginal contribution within the group. This contrasts with the ENSC rule, which distributes the non-separable component equally among all agents regardless of their individual cost structure.

An important interpretation of this rule is that it combines marginal responsibility with a proportional reward scheme based on individual cost potential. Agents whose individual cost is large relative to their separable cost receive a larger share of the non-separable component, reflecting a higher participation in the gains from cooperation. Using the decomposition of the non-separable cost, the SCRB rule can also be written as
\[
\begin{aligned}
SCRB_i(N,q,P)=&
a q_i
+
b\bigg(
q_N^{1-\eta}
-
q_{N\setminus\{i\}}^{1-\eta}
+
\frac{\left(
q_i^{1-\eta}
-
q_N^{1-\eta}
+
q_{N\setminus\{i\}}^{1-\eta}
\right)}{\sum_{j\in N}
\left(
q_j^{1-\eta}
-
q_N^{1-\eta}
+
q_{N\setminus\{j\}}^{1-\eta}
\right)
}
\bigg(q_N^{1-\eta}-\sum_{j\in N}(q_N^{1-\eta}-q_{N\setminus \{j\}}^{1-\eta}  ) \bigg) \bigg) .
\end{aligned}\]

From an economic perspective, the SCRB method introduces a purely proportional sharing of the interaction term, where weights are endogenously determined by the difference between individual and separable costs.

The allocation rules introduced above provide analytical expressions for distributing the total cost generated by group purchasing. To complement these theoretical results, the following example illustrates how the different
components of the cost decomposition and the corresponding allocation rules can be computed for a particular group-buying situation. In addition, it shows how the Equal Price, ENSC, and SCRB rules can be implemented computationally for the Schotanus price function. This computational perspective will also be useful in the following section, where we move from the allocation of a given total cost to the analysis of endogenous pricing.

\begin{example}\label{ex:computational_implementation}
Consider a group-buying situation with three agents,
$N=\{1,2,3\}$, whose individual demands are $
q_1=2, q_2=3$, and $q_3=5$. Besides, we take the price function of Schotanus specified by  \[
P(q)=10+\frac{100}{\sqrt{q}}.
\]

The aggregate demand of the group is $q_N=q_1+q_2+q_3=10$. Consequently, the group price and total cost are $P(q_N)=P(10)=10+\frac{100}{\sqrt{10}}
= 41.623$ and $C^{GB}(N,q,P)
=
P(q_N)q_N
= 416.228$. The individual costs of each of each agent and the corresponding separable costs are given in Table \ref{tab:indcostandSC}.

\begin{table}[h!]
    \centering
    \begin{tabular}{c|ccc}\hline
       & Agent 1 & Agent 2 & Agent 3 \\\hline\hline
       $C^{GB}(\{i\},q,P)$ & 161.421 & 203.205 & 273.607 \\
       $SC_i(N,q,P)$ & 53.385 & 81.653 & 142.621\\\hline
    \end{tabular}
    \caption{The individual costs and the separable costs.}
    \label{tab:indcostandSC}
\end{table}

Clearly, the individual costs  exceed the joint cost, since $638.233>416.228$. Besides, $\sum_{i\in N}SC_i(N,q,P)
= 277.659$, and the non-separable cost is 
\[
NSC(N,q,P)
=
C^{GB}(N,q,P)-\sum_{i\in N}SC_i(N,q,P)
= 138.569.\]

The Equal Price rule assigns to each agent a share of the total cost
proportional to her individual demand. Thus,
\[
EP(N,q,P)
=
P(q_N)q
=
(83.246,\;124.868,\;208.114).
\]

Under the ENSC rule, the non-separable cost is divided equally among the
three agents. Since $\frac{NSC(N,q,P)}{3}
= 76.601$, the resulting allocation is
\[
ENSC(N,q,P)=
(99.575,\;127.842,\;188.811).
\]

Finally, the SCRB rule distributes the non-separable cost according to the
relative gaps between individual costs and separable costs. These gaps are
approximately
\[
P(q_1)q_1-SC_1(N,q,P)= 108.036,
\]
\[
P(q_2)q_2-SC_2(N,q,P)= 121.552,
\]
and
\[
P(q_3)q_3-SC_3(N,q,P)= 130.986.
\]
Therefore, $SCRB(N,q,P)=(94.904,\;128.365,\;192.959)$.

Table \ref{tab:allocationsGB} summarizes the numerical results for the three allocation rules considered, using the cost-gap methodology and comparing the resulting cost allocations with the corresponding individual costs.

\begin{table}[h!]
    \centering
    \begin{tabular}{c|ccc|c}\hline
       $i$ & $EP_i$ & $ENSC_i$ & $SCRB_i$ & $C^{GB}(\{i\},q,P)$
       \\\hline\hline
       1 &83.246 & 99.575 & 94.904 & 161.421\\
       2 &  124.868 & 127.842 & 128.365 & 203.205\\
       3 &  208.114 & 188.811 & 192.959 & 273.607\\\hline
    \end{tabular}
    \caption{Cost allocation rules and individual costs.}
    \label{tab:allocationsGB}
\end{table}

As shown in the previous table, the allocations prescribed by the Equal Price rule, the ENSC-rule and the SCRB-rule yield cost shares that lie between the separable cost allocation, $SC(N,q,P)$, and the vector of individual costs, $P(q_i)q_i)_{i\in N}$, for every agent. The mentioned three allocation rules are also efficient, since
\[
\sum_{i\in N}EP_i(N,q,P)
=
\sum_{i\in N}ENSC_i(N,q,P)
=
\sum_{i\in N}SCRB_i(N,q,P)
=
C^{GB}(N,q,P)= 416.228.
\]
\noindent This example illustrates how the same total group-buying cost can lead to different individual payments depending on the principle used to distribute the non-separable component. The Equal Price rule preserves the group price for every unit purchased, the ENSC rule assigns the interaction component equally among agents, whereas the SCRB rule distributes this component according to the relative difference between individual and separable costs.
The comparison between these allocation principles provides the basis for the analysis that follows, where we move from the allocation of a given group cost to the determination of prices endogenously. \end{example}

\subsection{From allocation to endogenous pricing}

 A limitation of the standard pricing function $P(q)$ of \cite{schotanus2009unraveling} is that it depends solely on the individual quantity $q_i$, thereby implicitly assuming that each agent faces a price schedule independent of the composition of the purchasing group. However, in multi-agent group-buying environments, this assumption is overly restrictive. Prices are affected by aggregate demand and by the distribution of quantities across participants, since cooperation generates scale effects and redistributive externalities that depend on the entire demand profile $q=(q_i,q_{-i})$ (the amount of needing required by $i$ and agents in $N\setminus\{i\}$, respectively). In the remainder of the paper,  we formally denote $q_{-i}= (q_1, q_2,\cdots,q_{i-1},q_{i+1},\cdots, q_n)$ for all $i\in N$.

This motivates the introduction of \emph{variable discount  price functions} $\mathcal{P}:\mathbb{R}_+^{|N|}\longrightarrow \mathbb{R}$ of the form
\begin{equation}\label{generalpricef}
   \mathcal{P}(q_i,q_{-i})=a + b(q_i,q_{-i})\, q_i^{-\eta}, 
\end{equation}
where $\eta \in (0,1)$ preserves the elasticity structure of the original Schotanus price function, while $b$ is now allowed to depend on the full demand profile $(q_i,q_{-i})$. We refer to $b(q_i,q_{-i})$ as the variable discount factor. This approach  generalizes the standard unit price by allowing the price faced by agent $i$ to depend not only on her own demand $q_i$, but also on the demands of all other agents $q_{-i}$. Notice that the Equal Price rule can be interpreted as a variable discount price function because it can be rewritten as 
\begin{equation}
    \mathcal{P}(q_i,q_{-i})=a+b\bigg(\frac{q_i}{q_N}\bigg)^{\eta}\mbox{ } q_i^{-\eta}=P(q_N),
\end{equation}
with variable discount factor $b^{EP}(q_i,q_{-i})=b\bigg(\frac{q_i}{q_N}\bigg)^{\eta}>0$.

In what follows, we consider a general cost-allocation framework in which different allocation mechanisms in group-buying environments can be interpreted through their associated induced pricing representations. Specifically, each allocation rule $\chi$ generates a price function $P^{\chi}(q_i,q_{-i})$ that rationalizes the cost shares assigned to agents. This perspective provides a unifying framework for comparing different allocation principles in group-buying environments. 

\section{Bankruptcy-Driven Allocation Rules for Group-Buying}\label{sec3}

In this section, we analyze alternative ways of allocating the total cost  generated by cooperation among agents in $N$, adopting a bankruptcy perspective.  The key contribution of this framework is that cooperative purchasing endogenously generates a bankruptcy problem. Hence, classical bankruptcy rules can be used as cost-allocation mechanisms in group-buying environments.

A \emph{bankruptcy problem} (cf. \citealp{o1982problem}) involves a set of claimants $N$ whose total claims exceed the available estate. Formally, it is a triple $(E,d)$, where $E \in \mathbb{R}_+$ is the estate and $d \in \mathbb{R}^n$ is a vector of non-negative claims, with $0 \leq E \leq \sum_{j \in N} d_j$. The set of such problems is denoted by $B^N$.

For this purpose, we consider the bankruptcy perspective described below for any $(N,q,P)\in GB^N$, where the price function $P$ is specified as in (\ref{eq:psch}), following \cite{schotanus2009unraveling}.
\begin{definition}\label{defBANK}
Let $(N,q,P)\in GB^N$. The corresponding bankruptcy problem $(E^{GB},d^{GB})$ is defined by
 \[
 d^{GB}_i = C^{GB}(\{i\},q_i,P)-\big(C^{GB}(N,q,P)-C^{GB}(N\setminus\{i\},q_{-i},P)\big), \text{ for all } i \in N\]and\[
 E^{GB} = \sum_{i\in N}C^{GB}(\{i\},q_i,P)-C^{GB}(N,q,P).
 \]
\end{definition}
In our context, under the consideration of a Schotanus price function as the one in (\ref{eq:psch}), the vector of demands and the estate can be immediately written as
\begin{equation}\label{bankqi}
d^{GB}_i = b\bigg(q_i^{1-\eta} + q_{N\setminus\{i\}}^{1-\eta} - q_N^{1-\eta}\bigg), \text{for all } i \in N, \mbox{ and }
E^{GB} = b\bigg(\sum_{i\in N} q_i^{1-\eta} - q_N^{1-\eta}\bigg).
\end{equation}
Appendix \ref{appBANK} provides the verification that $(E^{GB},d^{GB})$ in this case constitutes a bankruptcy problem. {Moreover, it is worth noting that, for every $i\in N$, the claim $d_i^{GB}$ is non-decreasing in agent $i$'s demand $q_i$.}

Following the cost-gap methodology, the bankruptcy structure $(E^{GB},d^{GB})$ admits a natural economic interpretation. More specifically, $E^{GB}$ captures the economies of scale generated by cooperation and the resulting reduction in total expenditure. In this respect, it is structurally related to the non-separable cost component $NSC(N,q,P)$, since both reflect the cost differentials induced by coalition formation, scale effects, and the non-additivity of the underlying pricing structure. Likewise, for every $i\in N$, the demands of $i$ can be interpreted as agent $i$'s maximum claim on the total savings generated by cooperation: it represents the maximum amount that agent $i$ could reasonably claim from the cooperative surplus, whereas the bankruptcy rule determines the amount of the estate that the agent actually receives. Thus, both the estate and the individual claims are ultimately determined by the interaction between the separable and non-separable components induced by the group-buying technology. Assume, without loss of generality, that $d_1^{GB}\leq d_2^{GB}\leq\cdots\leq d_n^{GB}$.

The definition of procedures to allocate the estate among claimants is fundamental. A \emph{division rule} $f$ for a bankruptcy problem $(E,d)\in B^N$ is a function that assigns a vector $f(E,d)\in \mathbb{R}^n$ such that (i) $\sum_{i\in N} f_i(E,d)=E$, and (ii) $0 \leq f_i(E,d)\leq d_i$ for all $i\in N$. 
\textcolor{black}{Bankruptcy theory offers several classical rules for distributing a scarce estate among claimants with competing claims. Among them, the proportional rule (PROP) provides a natural benchmark by allocating the estate in proportion to the agents' claims (\citealp{o1982problem}). Other prominent rules embody different distributive principles. The constrained equal awards (CEA) rule seeks to equalize awards as much as possible, subject to the constraint that no claimant receives more than her claim, whereas the constrained equal losses (CEL) rule focuses on equalizing losses as much as possible, subject to non-negativity of the final allocations. The Talmud rule combines these two principles, applying an equal-awards logic when the estate is relatively small and an equal-losses logic when it is relatively large. Thus, CEA and CEL represent two dual egalitarian principles, while the Talmud rule provides an intermediate approach that switches between them according to the size of the estate relative to total claims (\citealp{aumann1985game}). Formally, these rules are defined as follows.}

\textcolor{black}{For any bankruptcy problem $(E,d)\in B^N$, the \emph{PROP rule} assigns to each claimant $i\in N$}

\begin{equation}\label{PROPrulebank}
    PROP_i(E,d)=\frac{d_i}{\sum_{i\in N}d_i}E.
\end{equation}\textcolor{black}{Hence, each claimant receives the same proportion of her claim, namely $\frac{E}{\sum_{j\in N}d_j}$}.

\textcolor{black}{We then consider the CEA, CEL, and Talmud rules, which depart from proportionality by applying alternative principles to the distribution of awards or losses.} The \emph{CEA rule} assigns to each claimant $i\in N$,
\begin{equation}\label{CEAorig}
    CEA_i(E,d)=\min \{d_i,\lambda\},
\end{equation}
where $\lambda$ satisfies $\sum_{i\in N}\min\{d_i,\lambda\}=E$. The \emph{CEL rule}  assigns to each claimant $i\in N$,
\begin{equation}\label{CELorig}
    CEL_i(E,d)=\max \{0,d_i-\lambda\},
\end{equation}
where $\lambda$ satisfies $\sum_{i\in N}\max\{0,d_i-\lambda\}=E$. {It turns out that}
\begin{equation}\label{CELorig2}
	CEL_i(E,d)=\max \{0,d_i-\lambda\}=d_i-CEA(\sum_{j\in N}d_j-E,d).
\end{equation} 
The \emph{TAL rule}, 
assigns to each claimant $i\in N$,
\[
TAL_i(N,d,E)=
\begin{cases}
CEA_i\!\left(E,\frac{d}{2}\right), & \text{if } E \leq \sum_{l\in N}\frac{d_l}{2}, \\[0.3em]
d_i - CEA_i\!\left(\sum_{i\in N}d_i-E,\frac{d}{2}\right), & \text{if } E \geq \sum_{l\in N}\frac{d_l}{2}.
\end{cases}
\]

For comprehensive surveys of bankruptcy problems and their associated division rules, including axiomatic characterizations based on natural and interpretable axioms, we refer the reader to \cite{HerreroVillar2001}, \cite{thomson2003axiomatic} and \cite{thomson2015axiomatic}.

Now, we consider a new perspective of rule for the class of group-buying multi-agent situations based on bankruptcy rules. Bankruptcy rules are not introduced as an alternative modelling assumption, but as a class of allocation mechanisms applied to the endogenous claim structure generated by the cooperative technology.

\begin{definition}\label{defpsiTU}
Let $(N,q,P)\in GB^N$ be a group-buying multi-agent situation  and take $\chi$ a bankruptcy rule. Thus, the $\chi$ rule in group-buying using bankruptcy is denoted by $\chi^{GB}=(\chi_i^{GB})_{i\in N}$ and is given  by 
\begin{equation}\label{propGBbank}
\chi^{GB}_i(N,q,P)=\chi_i(E^{GB},d^{GB}), \mbox{ for every }i \in N,
\end{equation}
and thus, the bankruptcy-driven cost  allocation for group buying is given by
\begin{equation} \label{BD-CA-GB}
\phi_i^{\chi,GB}(N,q,P)=C^{GB}(\{i\},q_i,P)- \chi_i(E^{GB},d^{GB}), \mbox{ for every }i \in N.
\end{equation}
\end{definition}

We should note that the bankruptcy-driven cost  allocation for group buying, $\phi_i^{\chi,GB}(N,q,P)$, is well-defined because $0 \leq \chi_i(E^{GB},d^{GB})\leq d_{i}^{GB} \mbox{ for every }i \in N$. Hence, $0 \leq SC_i(N,q,P)\leq \phi_i^{\chi,GB}(N,q,P) \leq C^{GB}(\{i\},q_i,P), \mbox{ for every }i \in N $ and $\sum_{i\in N}\phi_i^{\chi,GB}(N,q,P) =C^{GB}(N,q,P)$.
\medskip

In what follows, consider $(N,q,P)\in GB^N$ a group-buying multi-agent situation and the associated bankruptcy problem $(E^{GB},d^{GB})\in B^N$ specified in Definition \ref{defBANK}. The main purpose consists of identifying the variable discount price functions $\mathcal{P}^{\chi}(q_i,q_{-i})$ as the one in (\ref{generalpricef}) induced by each considered bankruptcy rule $\chi$ for any $(N,q,P)\in GB^N$. It is worth noting that the variable discount factor always reduces the initial discount factor; that is, $\mbox{ for every }i \in N,  b(q_i,q_{-i}) \leq b$ , by (\ref{BD-CA-GB}).   We refer to this function as the induced bankruptcy driven-group buying price function, henceforth, BD-GB price function. 


\subsection{On the proportional rule}
In this section, we propose a bankruptcy-driven allocation rule in group-buying based on the proportional rule for bankruptcy problems.

\begin{definition}
Let $(N,q,P)\in GB^N$ and consider the corresponding bankruptcy problem $(E^{GB},d^{GB})$. The proportional rule in group-buying via bankruptcy is denoted by $PROP^{GB}=(PROP_i^{GB})_{i\in N}$ and is given by 
\begin{equation}\label{propGBprop}
PROP^{GB}_i(N,q,P)=PROP_i(E^{GB},d^{GB}), \mbox{ for every }i \in N.
\end{equation}
\end{definition}

The proportional rule in (\ref{PROPrulebank})  assigns to each agent $i\in N$ the value
\begin{equation*}
PROP_i(E^{GB},d^{GB})
=b
\frac{q_i^{1-\eta}+q_{N\setminus\{i\}}^{1-\eta}-q_N^{1-\eta}}
{\sum_{j\in N}\big(q_j^{1-\eta}+q_{N\setminus\{j\}}^{1-\eta}-q_N^{1-\eta}\big)}
\left(\sum_{j\in N}q_j^{1-\eta}-q_N^{1-\eta}\right).
\end{equation*}

The associated bankruptcy-driven cost allocation $\phi_{i}^{PROP,GB}(N,q,P)$ is given by
\begin{equation}\label{phiPROP}
P(q_i)q_i- PROP_i(E^{GB},d^{GB})
= aq_i+bq_i^{1-\eta}
- b\,\frac{q_i^{1-\eta}+q_{N\setminus\{i\}}^{1-\eta}-q_N^{1-\eta}}
{\sum_{j\in N}\big(q_j^{1-\eta}+q_{N\setminus\{j\}}^{1-\eta}-q_N^{1-\eta}\big)}
\left(\sum_{j\in N}q_j^{1-\eta}-q_N^{1-\eta}\right).
\end{equation}

The BD-GB price function can be denoted as $\mathcal{P}^{PROP}$ for every $i\in N$. After some algebra, 
\[
\mathcal{P}^{PROP}(q_i,q_{-i})
=
a+b^{PROP}(q_i,q_{-i})\,q_i^{-\eta},
\]
where the variable discount factor is given by \[
b^{PROP}(q_i,q_{-i})
=
b\left(
1-
\frac{
q_i^{1-\eta}+q_{N\setminus\{i\}}^{1-\eta}-q_N^{1-\eta}
}{
\sum_{j\in N} (q_j^{1-\eta}
+
q_{N\setminus\{j\}}^{1-\eta}
-
q_N^{1-\eta})
}
\cdot
\frac{
\sum_{j\in N}q_j^{1-\eta}-q_N^{1-\eta}
}{
q_i^{1-\eta}
}
\right).
\]

Equivalently, using the bankruptcy notation, we obtain as variable discount factor
\[
0<b^{PROP}(q_i,q_{-i})
=
b\left(
1-
\frac{E^{GB}}{
\sum_{j\in N} d_j^{GB}
}
\cdot
\frac{
d_i^{GB}
}{
b\cdot q_i^{1-\eta}.
}
\right)
\]

The BD-GB price function associated with the proportional rule admits a natural economic interpretation. The parameter \(a\) represents the fixed baseline component of the unit price, which remains unaffected by cooperation. The term $b^{PROP}(q_i,q_{-i})\,q_i^{-\eta}$ captures the variable component of the price and incorporates the cooperative adjustment generated by the group-buying process. Furthermore, the original discount factor $b$ is scaled down by  $$\left(1-\frac{E^{GB}}{\sum_{j\in N} d_j^{GB}}\cdot\frac{d_i^{GB}}{b\cdot q_i^{1-\eta}} \right)$$ since that $$0<\frac{E^{GB}}{\sum_{j\in N} d_j^{GB}}<1\mbox{ and }0<\frac{d_i^{GB}}{b\cdot q_i^{1-\eta}}<1.$$

More precisely, the proportional bankruptcy rule modifies the elasticity-related coefficient \(b\) through an endogenous correction factor that depends on the agent's bankruptcy claim \(d_i^{GB}\). 
Cooperation operates by reducing the effective variable-price component faced by agents, while preserving the decreasing-demand structure of the original Schotanus pricing function.

\subsection{On the CEA rule}\label{sec:CEA}
Now, we propose a bankruptcy-driven allocation rule in group-buying based on the CEA rule for bankruptcy problems.

\begin{definition}
Let $(N,q,P)\in GB^N$ and consider the corresponding bankruptcy problem $(E^{GB},d^{GB})$. The constrained equal awards rule in group-buying via bankruptcy is denoted by $CEA^{GB}=(CEA_i^{GB})_{i\in N}$ and is given, for every $i\in N$, by 
\begin{equation}\label{propGBCEA}
CEA^{GB}_i(N,q,P)=CEA_i(E^{GB},d^{GB}).
\end{equation}
\end{definition}

We now derive the explicit expression of this rule and the corresponding cost allocation, distinguishing two cases. 

\begin{itemize}
\item[(a)] If $d^{GB}_1 \geq \frac{E^{GB}}{n}$, it holds that
\begin{equation*}
CEA_i(E^{GB},d^{GB}) = \frac{E^{GB}}{n}
= b\frac{\sum_{i\in N}q_i^{1-\eta}-q_N^{1-\eta}}{n}, \quad \forall i\in N.
\end{equation*}

Thus, the associated cost allocation $\phi^{CEA,GB}(N,q,P)$ is
\begin{equation*}
P(q_i)q_i - CEA_i(E^{GB},d^{GB})
= aq_i + b\cdot  q_i^{1-\eta}
- b\cdot \frac{\sum_{j\in N}q_j^{1-\eta}-q_N^{1-\eta}}{n},
\end{equation*}
that dividing by $q_i$, we obtain the following BD-GB price function:
\begin{align*}
    \mathcal{P}^{CEA}(q_i,q_{-i}) &
={a}+{b}^{CEA}(q_i,q_{-i})\frac{1}{q_i^{\eta}},
\end{align*}
where the variable discount factor is specified by $$b^{CEA}(q_i,q_{-i})=b\cdot \bigg(1-\frac{\sum_{j\in N}q_j^{1-\eta}-q_N^{1-\eta}}{n q_i^{1-\eta}}\bigg).$$ The term $1-\frac{\sum_{j\in N}q_j^{1-\eta}-q_N^{1-\eta}}{nq_i^{1-\eta}}$ is a discount factor that is less than one since that, for each $i\in N$,  $0<\sum_{j\in N}q_j^{1-\eta}-q_N^{1-\eta}<n q_i^{1-\eta}$.  Therefore, $0<b^{CEA}(q_i,q_{-i}).$

\item[(b)] Otherwise, let $k\in \{1,\ldots,n-1\}$ be such that
\begin{equation}\label{CEAb}
d^{GB}_k \leq \frac{b}{n-k}\left(\sum_{i\in N} q_i^{1-\eta} - q_N^{1-\eta} - \sum_{j=1}^k \bigg(q_j^{1-\eta} + q_{N\setminus\{j\}}^{1-\eta} - q_N^{1-\eta}\bigg)\right) < d^{GB}_{k+1}.
\end{equation}

The CEA allocation is given by
\begin{equation*}
CEA_i(E^{GB},d^{GB})=
\begin{cases}
b\bigg(q_i^{1-\eta} + q_{N\setminus\{i\}}^{1-\eta} - q_N^{1-\eta}\bigg), & \text{if } i \leq k,\\[6pt]
\displaystyle b\cdot \frac{\sum_{i\in N} q_i^{1-\eta} - q_N^{1-\eta} - \sum_{j=1}^k (q_j^{1-\eta} + q_{N\setminus\{j\}}^{1-\eta} - q_N^{1-\eta})}{n-k}, & \text{if } i>k.
\end{cases}
\end{equation*}

The associated cost allocation $\phi^{CEA,GB}(N,q,P)$ is obtained below. We distinguish two cases.

\begin{itemize}
\item[(i)] For $i \leq k$, the cost of agent i coincide with its separable cost specified in (\ref{SCrule}). Indeed,
\begin{align*}
P(q_i)q_i -  CEA_i(E^{GB},d^{GB})
&= aq_i + b q_i^{1-\eta} - b\big(q_i^{1-\eta}+q_{N\setminus\{i\}}^{1-\eta}-q_N^{1-\eta}\big) \\
&= aq_i + b\big(q_N^{1-\eta}-q_{N\setminus\{i\}}^{1-\eta}\big) \\
&= C^{GB}(N,q,P) - C^{GB}(N\setminus\{i\},q_{-i},P)\\
& = SC_i(N,q,P).
\end{align*}
This induces the following BD-GB price function $\mathcal{P}^{CEA}(q_i,q_{-i})=a+b^{CEA}(q_i,q_{-i})\frac{1}{q_i^{\eta}}$ such that the variable discount factor is
 $$0<b^{CEA}(q_i,q_{-i})=b\left(\frac{q_N^{1-\eta}-q_{N\setminus\{i\}}^{1-\eta}}{q_i^{1-\eta}}\right),$$
 by the increasing of the function $f(x)=x^{1-\eta}$. The variable discount factor reflects here the marginal contribution of agent $i$ to the efficiency gains of the grand coalition: $q_N^{1-\eta}-q_{N\setminus\{i\}}^{1-\eta}$ measures the additional collective efficiency generated by $i$, normalized by $q_i^{1-\eta}$. 

\item[(ii)] For $i > k$, $\phi^{CEA,GB}(N,q,P)$ is specified by
\begin{align*} 
P(q_i)q_i -  CEA_i(E^{GB},d^{GB})
&= aq_i + b q_i^{1-\eta}
- \frac{b}{n-k}\left( \sum_{j=k+1}^n q_j^{1-\eta}
+ \sum_{j=1}^k (q_N^{1-\eta}-q_{N\setminus\{j\}}^{1-\eta})
- q_N^{1-\eta} \right).
\end{align*}In this case, the BD-GB price function is $\mathcal{P}^{CEA}(q_i,q_{-i})=a+b^{CEA}(q_i,q_{-i})\frac{1}{q_i^{\eta}}$ with the variable discount factor   $$0<b^{CEA}(q_i,q_{-i})={b}\Bigg( 1-\frac{1}{q_i^{1-\eta}(n-k)}\left(
		\sum_{j=k+1}^n q_j^{1-\eta}
		+\sum_{j=1}^k\left(q_N^{1-\eta}-q_{N\setminus\{j\}}^{1-\eta}\right)
		-q_N^{1-\eta}\right)
		\Bigg),$$
        by the condition of $i>k$ and the inequality in (\ref{CEAb}). 
\end{itemize}

\end{itemize}

\subsection{On the CEL rule}\label{sec:CEL}

In this section, we propose a bankruptcy-driven allocation rule in group-buying based on the CEL rule for bankruptcy problems.

\begin{definition}
Let $(N,q,P)\in GB^N$ and consider the corresponding bankruptcy problem $(E^{GB},d^{GB})$. The constrained equal losses rule in group-buying via bankruptcy is denoted by $CEL^{GB}=(CEL_i^{GB})_{i\in N}$ and is given  by 
\begin{equation}\label{propGBcel}
CEL^{GB}_i(N,q,P)=CEL_i(E^{GB},d^{GB}), \mbox{ for every }i \in N.
\end{equation}
\end{definition}

The obtaining of the explicit expression of this rule follows a similar structure to the case of the CEA rule. The details, for the sake of readability, are included in Appendix \ref{CELbank}. We now derive the associated expression of price functions induced by the CEL rule, distinguishing two cases. 

From Equation~(\ref{CELorig2}) and considering the allocation prescribed by $CEA(	\sum_{i\in N}d^{GB}_i-E^{GB},d^{GB})$, we distinguish two scenarios:
\begin{itemize}
\item[(a)] If it holds $d_1^{GB} \geq \frac{\sum_{i\in N}d_i^{GB}-E^{GB}}{n}$, using Expression (\ref{CELaApp}) in Appendix B, the associated cost allocation $\phi^{CEL,GB}(N,q,P)$ is
\begin{equation}\label{casoAcel}
    \begin{aligned}
P(q_i)q_i -CEL_i(E^{GB},d^{GB}) &=aq_i+b\left(q_N^{1-\eta}-q_{N\setminus\{i\}}^{1-\eta}\right)+\frac{b}{n}\left(q_N^{1-\eta}+\sum_{j\in N}\left(q_{N\setminus\{j\}}^{1-\eta}-q_N^{1-\eta}\right)\right)\\
&=SC_i(N,q,P)+\frac{b}{n}\left(q_N^{1-\eta}+\sum_{j\in N}\left(q_{N\setminus\{j\}}^{1-\eta}-q_N^{1-\eta}\right)\right).
\end{aligned}
\end{equation}

We can assign a price function as
\[
\mathcal{P}^{CEL}(q_i,q_{-i})=a+b^{CEL}(q_i,q_{-i})\frac{1}{q_i^{\eta}}\]
where \[b^{CEL}(q_i,q_{-i})=\frac{b}{q_i^{1-\eta}}\left(q_N^{1-\eta}-q_{N\setminus\{i\}}^{1-\eta}+\frac{1}{n}\left(q_N^{1-\eta}+\sum_{j\in N}\left(q_{N\setminus\{j\}}^{1-\eta}-q_N^{1-\eta}\right)\right)\right),
\]{satisfying $0<b^{CEL}(q_i,q_{-i})$ because of the increasing of the function $f(x)=x^{1-\eta}$ and the fact that $0<\sum_{i\in N}d_i^{GB}-E^{GB}=q_N^{1-\eta}+\sum_{j\in N}\left(q_{N\setminus\{j\}}^{1-\eta}-q_N^{1-\eta}\right)$. }

The contribution $b^{CEL}(q_i,q_{-i})$ combines agent $i$'s marginal contribution to aggregate quality, $q_N^{1-\eta}-q_{N\setminus{i}}^{1-\eta}$, with a common adjustment term based on the aggregate quality $q_N^{1-\eta}$ and the marginal contributions of all agents, and normalizes the resulting expression by $q_i^{1-\eta}$.

\item[(b)] Otherwise, let $k\in \{1,\ldots,n-1\}$ be such that 
\begin{equation}\label{CELB}
d_k^{GB} \leq \frac{b}{n-k} \bigg(q_N^{1-\eta}+\sum_{j\in N} (q_{N\setminus\{j\}}^{1-\eta}-q_N^{1-\eta})-\sum_{j=1}^k (q_j^{1-\eta}+q_{N\setminus\{j\}}^{1-\eta}-q_N^{1-\eta})\bigg)< d_{k+1}^{GB},
\end{equation}

From the associated cost allocation given in (\ref{CELbApp}), the associated cost allocation $\phi^{CEL,GB}(N,q,P)$ and the BD-GB price functions can be considered when distinguishing by agent. Thus, 	
\begin{itemize}
\item \textit{Case $i\leq k$}. We have that $P(q_i)q_i - CEL_i(E^{GB},d^{GB}) = P(q_i)q_i$, and, hence, 
$\mathcal{P}^{CEL}(q_i,q_{-i})=a+b^{CEA}(q_i,q_{-i})\frac{1}{q_i^{\eta}}$ with $b^{CEL}(q_i,q_{-i}) = b$ because $CEL_i(E^{GB},d^{GB})=0$.
\item \textit{Case $i>k$}. We obtain that 
\begin{align*}
P(q_i)q_i -  CEL_i(E^{GB},d^{GB})
&= aq_i+bq_i^{1-\eta}-b\bigg(q_i^{1-\eta}+q_{N\setminus\{i\}}^{1-\eta}-q_N^{1-\eta}-	\\
&\hspace{3.5 cm}-\frac{q_N^{1-\eta}- \sum_{j=k+1}^n (q_N^{1-\eta}-q_{N\setminus\{j\}}^{1-\eta})-\sum_{j=1}^k q_j^{1-\eta}}{n-k}\bigg)\\
& =SC_i(N,q,P) +b\bigg(\frac{q_N^{1-\eta}- \sum_{j=k+1}^n (q_N^{1-\eta}-q_{N\setminus\{j\}}^{1-\eta})-\sum_{j=1}^k q_j^{1-\eta}}{n-k}\bigg).
\end{align*}

 and that $\mathcal{P}^{CEL}(q_i,q_{-i})=a+b^{CEL}(q_i,q_{-i})\frac{1}{q_i^{\eta}}$, with $$b^{CEL}(q_i,q_{-i})
=
b\bigg(\frac{q_N^{1-\eta}-q_{N\setminus\{i\}}^{1\eta}}{q_i^{1-\eta}}+\frac{q_N^{1-\eta}-\sum_{j=1}^k q_j^{1-\eta}+\sum_{j=k+1}^n \left(q_{N\setminus\{j\}}^{1-\eta}-q_{N\setminus\{i\}}^{1-\eta}\right)}{\left(n-k\right)q_i^{1-\eta}}\bigg),$$
\end{itemize}
satisfying that $0<b^{CEL}(q_i,q_{-i})$ by inequality (\ref{CELB}) and the increasing of the function $f(x)=x^{1-\eta}$.  The expression of $b^{CEL}(q_i,q_{-i})$ shows that the CEL rule changes the original coefficient $b$ according to the demands of all agents. In particular, $b^{CEL}$ depends on the total demand $q_N$, on the demand of agent $i$, and on the total demand that remains when some agents are excluded. Therefore, the adjustment faced by agent $i$ depends not only on its own demand, but also on how its demand compares with the demands of the other agents. The parameter $k$ determines how the first $k$ agents enter the expression, while the remaining $n-k$ agents are treated together through the denominator. The basic structure of the Schotanus pricing function is therefore preserved, but the coefficient becomes dependent on the distribution of demand across agents.
\end{itemize}

\subsection{On the Talmud rule}

Finally, we propose a bankruptcy-driven allocation rule in group-buying based on the Talmud rule for bankruptcy problems.

\begin{definition}
Let $(N,q,P)\in GB^N$ and consider the corresponding bankruptcy problem $(E^{GB},d^{GB})$. The Talmud rule in group-buying via bankruptcy is denoted by $TAL^{GB}=(TAL_i^{GB})_{i\in N}$ and is given by 
\begin{equation}\label{propGBtal}
TAL^{GB}_i(N,q,P)=TAL_i(E^{GB},d^{GB}), \mbox{ for every }i \in N.
\end{equation}
\end{definition}

The obtaining of the explicit expression of this rule for any $(N,q,P)\in GB^N$ follows a similar structure to the previous cases. The complete version of results, for the sake of readability, are included in Appendix \ref{TalmudBank}. This allows the determination of the associated price function induced by the Talmud rule. We distinguish two scenarios:

\begin{itemize}
\item[(a)] If $E^{GB} \leq \frac{1}{2}\sum_{i\in N} d_i^{GB}$, take $k\in\{0,\ldots,n-1\}$ with $d_0^{GB}=0$ such that
\begin{equation}\label{TalmudA}
\frac{d^{GB}_k}{2}
\leq \frac{b}{n-k}\left(\bigg(\sum_{i\in N} q_i^{1-\eta} - q_N^{1-\eta}\bigg)-\sum_{j=1}^k \frac{\bigg(q_j^{1-\eta} + q_{N\setminus\{j\}}^{1-\eta} - q_N^{1-\eta}\bigg)}{2}\right)
<\frac{d^{GB}_{k+1}}{2}.
\end{equation}

Using the expression of the Talmud allocation in Equation (\ref{TalaApp}) in Appendix \ref{TalmudBank},
 $\phi^{TAL,GB}(N,q,P)$ admits two expressions according to the related position of agent $i$. 
\begin{itemize}
\item If $i\leq k$,
\begin{equation*}
P(q_i)q_i-TAL_i(E^{GB},d^{GB})
= aq_i+bq_i^{1-\eta}
-\frac{b}{2}\left(q_i^{1-\eta}+q_{N\setminus\{i\}}^{1-\eta}-q_N^{1-\eta}\right).
\end{equation*}

\item If $i>k$,
\begin{equation*}
P(q_i)q_i-TAL_i(E^{GB},d^{GB})
= aq_i+bq_i^{1-\eta}
-\frac{b}{n-k}\left( \sum_{j\in N} q_j^{1-\eta} - q_N^{1-\eta}-\frac{1}{2}\sum_{j=1}^k \bigg(q_j^{1-\eta} + q_{N\setminus\{j\}}^{1-\eta} - q_N^{1-\eta}\bigg)\right).
\end{equation*}
\end{itemize}

Thus, the associated BD-GB price function takes the form $\mathcal{P}^{TAL}(q_i,q_{-i})
=a
+ b^{TAL}(q_i,q_{-i}) q_i^{-\eta}$ such that:
\begin{itemize}
\item[(i)] For $i \leq k$, it holds, by the increasing of the function $f(x)=x^{1-\eta}$, that $$0<b^{TAL}(q_i,q_{-i})
=
\frac{b}{2}\left(1+\frac{q_N^{1-\eta}-
q_{N\setminus\{i\}}^{1-\eta}
}{
q_i^{1-\eta}
}
\right).$$

\item[(ii)] For $i > k$, we obtain, by (\ref{TalmudA}),
\[
0<b^{TAL}(q_i,q_{-i})
=
b\left(
1
-
\frac{
 \bigg(\sum_{j\in N} q_j^{1-\eta} - q_N^{1-\eta}\bigg)-\frac{1}{2}\sum_{j=1}^k \bigg(q_j^{1-\eta} + q_{N\setminus\{j\}}^{1-\eta} - q_N^{1-\eta}\bigg)
}{
(n-k)\,q_i^{1-\eta}}
\right).
\]
\end{itemize}

\item[(b)] If $E^{GB} \geq \frac{1}{2}\sum_{i\in N} d_i^{GB}$, take $k\in \{0,\cdots,n-1\}$ with $d_0^{GB}=0$ such that:

\begin{equation}\label{condTALb}
\frac{d^{GB}_k}{2}
\leq \frac{b}{n-k}\left(\sum_{j\in N}\bigg( q_{N\setminus\{j\}}^{1-\eta} - q_N^{1-\eta}\bigg) + q_N^{1-\eta} - \sum_{j=1}^k \frac{\bigg(q_j^{1-\eta} + q_{N\setminus\{j\}}^{1-\eta} - q_N^{1-\eta}\bigg)}{2}\right)
<\frac{d^{GB}_{k+1}}{2}.
\end{equation}
Then, for every $i\in N$, the corresponding bankruptcy-driven cost allocation $\phi^{TAL,GB}(N,q,P)$  is given, for every $i\in N$,  by the following expressions, that depend on the relative position of $i$.
\begin{itemize}
    \item If $i\leq k$, it is given by
\begin{equation*}
P(q_i)q_i-TAL_i(E^{GB},d^{GB})
= aq_i+bq_i^{1-\eta}
-\frac{b}{2}\left(q_i^{1-\eta}+q_{N\setminus\{i\}}^{1-\eta}-q_N^{1-\eta}\right).
\end{equation*}
\item Otherwise, if $i>k$,
\hspace*{1 cm}\begin{align*}
P(q_i)q_i-TAL_i(E^{GB},d^{GB})
&= aq_i+bq_i^{1-\eta}-b (q_i^{1-\eta} + q_{N\setminus\{i\}}^{1-\eta} - q_N^{1-\eta}\big) +\\
&\hspace{0.25 cm}+\frac{b}{n-k}\left( q_N^{1-\eta}+\sum_{j\in N}\bigg( q_{N\setminus\{j\}}^{1-\eta} - q_N^{1-\eta}\bigg)  - \sum_{j=1}^k \frac{\bigg(q_j^{1-\eta} + q_{N\setminus\{j\}}^{1-\eta} - q_N^{1-\eta}\bigg)}{2}\right)\\
&= C^{GB}(N,q,P)-C^{GB}(N\setminus\{i\},q_{-i},P) +\\
&\hspace{0.25 cm}+\frac{b}{n-k}\left( q_N^{1-\eta}+\sum_{j\in N}\bigg( q_{N\setminus\{j\}}^{1-\eta} - q_N^{1-\eta}\bigg)  - \sum_{j=1}^k \frac{\bigg(q_j^{1-\eta} + q_{N\setminus\{j\}}^{1-\eta} - q_N^{1-\eta}\bigg)}{2}\right).
\end{align*}

Thus, the associated BD-GB price function $\mathcal{P}^{TAL}(q_i,q_{-i})
=a
+ b^{TAL}(q_i,q_{-i}) q_i^{-\eta}$ is given as follows.
\begin{itemize}
\item[(i)] For $i \leq k$, we have that \[0<b^{TAL}(q_i,q_{-i})
=
\frac{b}{2}\left(1+\frac{
q_N^{1-\eta}-q_{N\setminus\{i\}}^{1-\eta}}{
q_i^{1-\eta}
}
\right),\]

once again by the increasing of function $f(x)=x^{1-\eta}$.

\item[(ii)] For $i > k$, we obtain, by (\ref{condTALb}), that
%
\begin{align*}
0<b^{TAL}(q_i,q_{-i})
&=
b\Bigg(
\frac{q_N^{1-\eta}- q_{N\setminus\{i\}}^{1-\eta}
}{
q_i^{1-\eta}
}
\\
& \quad +
\frac{
\displaystyle
 q_N^{1-\eta}+\sum_{j\in N}\bigg( q_{N\setminus\{j\}}^{1-\eta} - q_N^{1-\eta}\bigg)
}{
(n-k)\,q_i^{1-\eta}
}-
\frac{
\displaystyle
\sum_{j=1}^{k} {\bigg(q_j^{1-\eta} + q_{N\setminus\{j\}}^{1-\eta} - q_N^{1-\eta}\bigg)}
}{
2(n-k)\,q_i^{1-\eta}
}
\Bigg).
\end{align*}

\end{itemize}
\end{itemize}\end{itemize}
In both cases, the Talmud rule preserves the baseline Schotanus component although all redistributive effects are captured exclusively through \(b^{TAL}(q_i,q_{-i})\). Under low estate, agents \(i \leq k\) receive reductions proportional to half of their claims and marginal contributions, whereas agents $i>k$ are subject to an adjustment determined by the remaining estate and scaled by demand, reflecting partial compensation under scarcity. Under high estate, the same structure applies, but the adjustment for agents $i>k$ incorporates the residual imbalance after the partial allocation to the first $k$ agents. Hence, the rule modifies only the demand-dependent component of prices, leaving the fixed term unchanged and generating regime-dependent elastic adjustments.

\subsection{Comparison and relationships among allocation rules}\label{sec35}

The bankruptcy formulation provides an alternative perspective on cost allocation in group-buying situations. Once the cooperative surplus generated by group purchasing is represented as a bankruptcy estate, different bankruptcy principles give rise to different cost-sharing rules. This section compares these bankruptcy-driven allocations with the cost-gap-based rules introduced above and investigates the relationships between the two approaches.

We begin with a numerical example that illustrates the different allocations generated by the two approaches. The example also serves to identify potential correspondences between particular bankruptcy rules and cost-gap-based allocation rules. These correspondences are subsequently examined analytically.

\begin{example}\label{ex:bankruptcy_implementation}
We revisit the group-buying situation considered in Example \ref{ex:computational_implementation}. The price function is given by the Schotanus specification $P(q)=10+\frac{100}{\sqrt{q}}$. Recall that $C^{GB}(N,q,P)
= 416.228$.

The associated bankruptcy problem is characterized by an estate of $E^{GB}=222.006$ and claims $d^{GB}$, as reported in Table \ref{tab:data3.7}. Such table also reports the demands and the individual and the separable costs.

\begin{table}[h!]
    \centering
    \begin{tabular}{c|cc|cc}\hline
       $i$ & $q_i$  & $d_i^{GB}$ & $C^{GB}(\{i\},q,P)$ &$SC_i$\\\hline\hline
       1    &    2& 108.0363 &    161.4214& 53.385\\
       2& 3 & 121.5524 & 203.2051 & 81.653\\
       3&5&130.9858&273.6068& 142.621\\\hline
    \end{tabular}
    \caption{Demands, claims, individual and separable costs per agent.}
    \label{tab:data3.7}
\end{table}

Below, we display the Proportional rule, the CEA rule, the CEL rule and the Talmud rule for group-buying in Table \ref{tab:allocationsGB2}. We also show the cost-gap-based allocations described in Section~\ref{costgap}. A complete comparison of the two approaches is also reported.
\begin{table}[h!]
    \centering
    \begin{tabular}{c|ccc|ccc}\hline
   &  \multicolumn{3}{c}{Cost gap-based allocations}&\multicolumn{3}{|c}{Bankruptcy-driven  allocations}\\
    \hline
       $i$   & $EP_i$ & $ENSC_i$ &$SCRB_i$  & $\phi_i^{PROP,GB}$ & $\phi_i^{CEA,GB}$ & $\phi_i^{CEL,GB}=\phi_i^{TAL,GB}$
\\\hline\hline
       1  & 83.246 & 99.575 & 94.904 & 94.904 & 87.420 & 99.575\\
       2   & 124.868 & 127.842 & 128.365 & 128.365 & 129.203 & 127.842\\
       3 & 208.114 & 188.811 & 192.959   & 192.959 &199.605 & 188.811\\\hline
    \end{tabular}
    \caption{Comparison of different cost allocation proposals.}
    \label{tab:allocationsGB2}
\end{table}

Several observations emerge from this comparison. First,  bankruptcy-driven rules usually generate different distributions of the total cost. In particular, CEA assigns a relatively lower cost to agent $1$ and a relatively higher cost to agent $3$ than the other bankruptcy-driven rules considered in the example. Second, some bankruptcy-driven allocations coincide with cost-gap-based allocations in this particular group-buying situation. The proportional rule yields the same allocation as SCRB, while CEL, Talmud, and ENSC also happen to yield the same allocation. These numerical correspondences motivate the theoretical analysis that follows, where we determine which relationships hold generally and which depend on the particular group-buying situation.\end{example}

The numerical example above shows that, in some group-buying situations, bankruptcy-driven cost allocation rules based on the $CEA$, $CEL$, and Talmud rules may yield the same allocation. The following example demonstrates that this need not always be the case by considering a situation in which the three rules generate different allocations.
\begin{example}
Consider the group-buying situation with three agents, $N=\{1,2,3\}$, whose individual demands are given in the second column of Table \ref{tab:data3.11}. Besides, we take the price function of Schotanus specified by $P(q)=1.43+\frac{11.48}{q^{0.01}}$. In this case, we obtain $C^{GB}(N,q,P) = 251.424$ and $NSC(N,q,P)=1.708$.

The corresponding bankruptcy problem $(E^{GB},d^{GB})$ is characterized by an estate of $E^{GB}=2.238$ and the claims $d^{GB}$ given below. The table also presents the individual and separable costs, together with the resulting bankruptcy-based cost allocations under consideration.

\begin{table}[h!]
\centering
\begin{tabular}{c|cccc|ccc}\hline
$i$ & $q_i$ & $d_i^{GB}$ & $C^{GB}(\{i\},q,P)$ & $SC_i$ &
$\phi_i^{CEA,GB}$ & $\phi_i^{CEL,GB}$ & $\phi_i^{TAL,GB}$
\\\hline\hline
1 & 7  & 1.448 & 88.821 & 87.373 & 88.075 & 87.942 & 87.990\\
2 & 10 & 1.550 & 126.487 & 124.937& 125.741 & 125.506 & 125.554\\
3 & 3  & 0.948 & 38.354 &  37.406& 37.608 & 37.975 & 37.880\\\hline
\end{tabular}
\caption{Demands, claims, individual and separable costs,  and BD-GB rules}
\label{tab:data3.11}
\end{table}
It follows that the three bankruptcy-driven cost allocation rules generate different distributions of the cooperative savings among the agents.\end{example}


The numerical example motivates the following question: to what extent can cost-gap-based allocation rules be represented as bankruptcy-driven cost allocations? We first establish that the correspondence observed between the proportional bankruptcy rule and SCRB is not specific to the numerical examples, but holds for every group-buying situation in $GB^N$.

 \begin{proposition}
  Let $(N,q,P)\in GB^N$.   The  bankruptcy-driven cost allocation $\phi^{PROP,GB}(N,q,P)$ corresponds to the $SCRB(N,q,P)$, i.e. $\phi^{PROP,GB}(N,q,P)=SCRB(N,q,P).$
 \end{proposition}
 \begin{proof}
    Take $(N,q,P)\in GB^N$. For every $i\in N$, by Expression in (\ref{phiPROP}), $\phi^{PROP,GB}_i(N,q,P)$ can be formalized as follows.
	\begin{equation*}
		\phi^{PROP,GB}_i(N,q,P)=aq_i+bq_i^{1-\eta}-b\frac{q_i^{1-\eta}+q_{N\setminus\{i\}}^{1-\eta}-q_N^{1-\eta}}{\sum_{j\in N}(q_j^{1-\eta}+q_{N\setminus\{j\}}^{1-\eta}-q_N^{1-\eta})}\Big(\sum_{j\in N}q_j^{1-\eta}-q_N^{1-\eta}\Big).
	\end{equation*}
	On the other hand, using Expression (\ref{SCRBrule}), we have that
	\begin{equation*}
SCRB_i(N,q,p)=aq_i+b\Big( q_N^{1-\eta}-q_{N\setminus\{i\}}^{1-\eta}+\frac{q_i^{1-\eta}+q_{N\setminus\{i\}}^{1-\eta}-q_N^{1-\eta}}{\sum_{j\in N}(q_j^{1-\eta}+q_{N\setminus\{j\}}^{1-\eta}-q_N^{1-\eta})}\Big(q_N^{1-\eta}-\sum_{j\in N}(q_N^{1-\eta}-q_{N\setminus\{j\}}^{1-\eta})\Big)\Big),
	\end{equation*}
    for every $i\in N$. Then, subtracting the two expressions and simplifying, we have
	\begin{align*}
			\phi^{PROP,GB}_i(N,q,P)-SCRB_i(N,q,p)&=b\Big( q_i^{1-\eta}+q_{N\setminus\{i\}}^{1-\eta}-q_N^{1-\eta}\Big)-\\
			&\hspace{1 cm}-b\frac{q_i^{1-\eta}+q_{N\setminus\{i\}}^{1-\eta}-q_N^{1-\eta}}{\sum_{j\in N}(q_j^{1-\eta}+q_{N\setminus\{j\}}^{1-\eta}-q_N^{1-\eta})}\sum_{j\in N}(q_j^{1-\eta}+q_{N\setminus\{j\}}^{1-\eta}-q_N^{1-\eta})\\
			&\hspace{1 cm}=0,
	\end{align*}
and therefore $\phi^{PROP,GB}(N,q,P)=SCRB(N,q,P)$.
 \end{proof}

We next turn to the comparison between the $\phi^{CEA,GB}$ and ENSC allocations. Unlike the proportional rule, the relationship between these two allocations is not characterized by a general equivalence. Instead, their relative values depend on the structure of the bankruptcy problem and, in particular, on the position of the claims relative to the equal-awards level under the CEA rule.

Using the same Cases as in Section \ref{sec:CEA}, we first consider Case (a), in which $d_1^{GB}\geq \frac{E^{GB}}{n}$ 
\begin{align*}
	\phi^{CEA,GB}_i(N,q,P)-	ENSC_i(N,q,p)&=b\Big( q_i^{1-\eta}+q_{N\setminus\{i\}}^{1-\eta}-q_N^{1-\eta}\Big)-\frac{b}{n}\sum_{j\in N}(q_j^{1-\eta}+q_{N\setminus\{j\}}^{1-\eta}-q_N^{1-\eta})\\
	&=d^{GB}_i-\frac{1}{n}\sum_{j\in N}d^{GB}_j
\end{align*}	
	This difference can be positive or negative depending on the relationship between $d_i$ and the average of all values in $d$. For Case (b), when for every $i\in N$, we take
 $k\in \{1,\ldots,n-1\}$  such that
\begin{equation}
d^{GB}_k \leq \frac{b}{n-k}\left(\bigg(\sum_{i\in N} q_i^{1-\eta} - q_N^{1-\eta}\bigg) - \sum_{j=1}^k \bigg(q_j^{1-\eta} + q_{N\setminus\{j\}}^{1-\eta} - q_N^{1-\eta}\bigg)\right) < d^{GB}_{k+1},
\end{equation}
	we conclude the following statements.
\begin{itemize}
    \item If $i\leq k$, it readily follows that $\phi^{CEA,GB}_i(N,q,p)=SC_i(N,q,p)<	ENSC_i(N,q,p)$.
    \item  If $i> k$, we do not obtain a conclusive relation between the allocations $\phi^{CEA,GB}_i(N,q,p)$ and $ENSC_i(N,q,p)$.    The next example illustrates this point.
\end{itemize}

\begin{example}Consider the group-buying situation with seven agents, $N=\{1,2,3,4,5,6,7\}$, whose individual demands are given in the second column of Table \ref{tab:data3.9}. Besides, we take the price function of Schotanus specified by  $P(q)=10+\frac{100}{\sqrt{q}}$. Now, we obtain $C^{GB}(N,q,P) = 4395.295$ and $NSC(N,q,P)=774.341$.

The corresponding bankruptcy problem $(E^{GB},d^{GB})$ is characterized by an estate of $E^{GB}=2302.193$ and the claims $d^{GB}$ given below. We also present the corresponding allocations obtained using the CEA rule, along with additional numerical results.

\begin{table}[h!]
    \centering
    \begin{tabular}{c|cccc|ccc}\hline
       $i$ & $q_i$  & $d_i^{GB}$ & $C^{GB}(\{i\},q,P)$ & $SC_i$ & $CEA_i$&$\phi_i^{CEA,GB}$ & $ENSC_i$       \\\hline\hline
 1    &   2 &135.369 & 161.421&26.052& 135.369&26.052& 136.672\\
 2&8& 258.499&        362.843&104.344 & 258.499&  104.344 & 214.964\\
3&20& 385.657& 647.214& 261.557& 381.665& 265.549 &  372.177\\
4& 38 &497.376&996.441& 499.065& 381.665& 614.777&  609.686\\
5&48& 540.855&1172.820&631.965 &  381.665&791.155& 742.585\\
6&74&619.152&1600.233&981.081& 381.665& 1218.568&  1091.701\\
7&84&639.626&1756.515& 1116.890& 381.665&1374.850& 1227.510\\\hline
    \end{tabular}
    \caption{Comparison of  $\phi^{CEA,GB}$ and $ENSC$ rules and more numerical results.}
    \label{tab:data3.9}
\end{table}

	On the other hand, we compute the ENSC-rule in (\ref{ENSCrule}) for $(N,q,P)$. The comparison illustrates that the difference between $\phi^{CEA,GB}$ and ENSC is not uniform across agents. In particular, agents $1$, $2$, and $3$ pay less under the $\phi^{CEA,GB}$ allocation than under ENSC. By contrast, agents $4,5,6$ and $7$ pay more under $\phi^{CEA,GB}$ than under ENSC. Thus, even within the same group-buying situation, the two rules redistribute the cooperative gains differently across agents. This example illustrates why, in Case (b), no general ordering between the $\phi^{CEA,GB}$ and ENSC cost shares can be established for the agents whose claims are not binding.\end{example}

Finally, we next compare the $\phi^{CEL,GB}$ and ENSC allocations. From equations  (\ref{ENSCrule0}) and (\ref{casoAcel}), we observe that the two rules coincide, so that $\phi^{CEL,GB}(N,q,P)=ENSC(N,q,P)$. This coincidence, however, does not extend to Case (b). In this case, the two allocation procedures may distribute the cooperative gains differently across agents, as illustrated by the following example.

\begin{example}Consider the group-buying situation with seven agents,
$N=\{1,2,3,4,5,6,7\}$, whose individual demands are reported in the second column of Table \ref{tab:data3.10}. We consider the price function proposed by Schotanus  $P(q)=10+\frac{100}{{q}^{0.1}}$. In this case,  $C^{GB}(N,q,P) = 35991.750$ and $NSC(N,q,P)=242.024$.

\begin{table}[h!]
    \centering
    \begin{tabular}{c|cccc|ccc}\hline
       $i$ & $q_i$  & $d_i^{GB}$ & $C^{GB}(\{i\},q,P)$ & $SC_i$  & $CEL_i$&$\phi_i^{CEL,GB}$ & $ENSC_i$       \\\hline\hline
1   &     2  & 86.234     &    206.607 & 206.767 & 0.000 & 206.607 &  364.611\\
2   &    10  &291.874     &    894.328 & 1033.943 & 16.034 & 878.294 &   846.693\\
3   &    11  &312.692     &    975.473 & 1137.353& 36.852 & 938.621 &907.020\\
4   &    32  &649.529     &   2582.742 & 3309.613  & 373.689  &2209.053&  2177.452\\
5   &    90 &1157.836     &   6638.763 &9316.336& 881.996  &5756.766& 5725.165\\
6   &   200 &1356.697     &  13774.080 & 20745.714& 1080.857&12693.223& 12661.622\\\hline
    \end{tabular}
    \caption{Comparison of  $\phi^{CEL,GB}$ and $ENSC$ rules and more numerical results.}
    \label{tab:data3.10}
\end{table}

The associated bankruptcy problem has estate
$E^{GB}=2389.428$ and claims $d^{GB}$ as specified in Table~\ref{tab:data3.10}. The corresponding constrained equal losses allocation, $CEL(E^{GB},d^{GB})$, and the induced allocation $\phi^{CEL,GB}$ are reported in the table, together with the $ENSC(N,q,P)$ allocation.

This example shows that, once some claims are binding, $\phi^{CEL,GB}$ and ENSC no longer coincide. In particular, buyer $1$ pays less under the CEL allocation than under ENSC (that is even greater than its individual cost), whereas buyers $2,\ldots,6$ pay less under ENSC than under CEL. Thus, the difference between the two rules is not merely a consequence of the common total cost, but reflects a different distribution of the gains from group purchasing across agents. Furthermore, unlike ENSC, $\phi^{CEL,GB}$ does not necessarily assign to every agent $i\in N$ a cost greater than their individual cost $C^{GB}(\{i\},q,P)$.
\end{example}


\section{Numerical and simulation analysis}\label{sec4}

This section presents a stylized numerical example based on simulated data, allowing us to explore the behavior and distributional implications of the proposed allocation mechanisms.

\textcolor{black}{%
We provide an economically meaningful interpretation of the numerical analysis, by considering a purchasing group formed by 100 independent small and medium-sized distributors of electrical supplies. All distributors periodically purchase the same standardized electrical component (e.g., a circuit breaker) from a common supplier, but their individual requirements differ substantially because of differences in firm size, customer base, and local market demand. If each distributor purchases independently, the unit price is determined by its own order quantity. Alternatively, the distributors may pool their orders through the purchasing group, allowing the aggregate demand to qualify for a more favorable quantity discount. The resulting reduction in total procurement cost creates cooperative savings that must subsequently be distributed among the participating buyers.}

\textcolor{black}{The numerical instance considered below is not intended to reproduce a particular purchasing agreement. Rather, it represents a stylized but realistic group-buying situation designed to capture the substantial demand heterogeneity typically observed among firms of different sizes. The simulated demands can therefore be interpreted as the quantities required by each one of the 100 distributors during a common purchasing period (for example, over a one-month or three-month period)\footnote{{The same setting may also represent centralized procurement within a single organization, where the buyers correspond to different branches or business units whose demands are aggregated to obtain quantity discounts.}}. In this setting, buyers with relatively small demands may be viewed as local distributors serving a limited customer base, whereas buyers with larger demands represent regional distributors supplying a broader network of customers.}


Consider a population of $n=100$ potential buyers interested in purchasing the same product from a supplier. Buyers differ in the quantity they demand, but they interact with the same group-buying scheme and, consequently, face the same underlying Schotanus price function. The purpose of this example is to illustrate how heterogeneity in individual demand can affect to the different bankruptcy-driven cost allocations analysed along this paper. 

First, we build the associated group-buying multi-agent situation $(N,q,P)$.  Fix $N=\{1,\cdots,100\}$. Let $q_i$ denote the demand of buyer $i$, for each $i\in N$. We generate a heterogeneous population by randomly selecting 100 distinct integer demand levels from the interval $[1000,3000]$. Hence, 
with $q_i\neq q_j$ for $i\neq j$. The resulting population represents buyers with substantially different purchasing requirements. In particular, the dispersion in $q_i$ allows the model to capture situations in which some buyers demand relatively small quantities, whereas others are interested in considerably larger purchases.

\footnotesize
\begin{table}[h!]
    \centering\resizebox{0.7\textwidth}{!}{
    \begin{tabular}{cc|cc|cc|cc|cc}\hline
     Buyer $i$ & $q_i$   & Buyer $i$ & $q_i$   & Buyer $i$ & $q_i$   & Buyer $i$ & $q_i$  & Buyer $i$ & $q_i$  \\\hline\hline
   1 & 1021 & 21 & 1470 & 41 & 1877 & 61 & 2327 & 81 & 2700\\
2 & 1036 & 22 & 1484 & 42 & 1888 & 62 & 2330 & 82 & 2707\\
3 & 1038 & 23 & 1498 & 43 & 1929 & 63 & 2349 & 83 & 2741\\
4 & 1039 & 24 & 1525 & 44 & 1930 & 64 & 2368 & 84 & 2748\\
5 & 1083 & 25 & 1531 & 45 & 1974 & 65 & 2394 & 85 & 2752\\
6 & 1110 & 26 & 1536 & 46 & 1982 & 66 & 2413 & 86 & 2765\\
7 & 1128 & 27 & 1553 & 47 & 1995 & 67 & 2425 & 87 & 2790\\
8 & 1192 & 28 & 1555 & 48 & 2016 & 68 & 2427 & 88 & 2798\\
9 & 1247 & 29 & 1581 & 49 & 2052 & 69 & 2445 & 89 & 2814\\
10 & 1251 & 30 & 1596 & 50 & 2068 & 70 & 2515 & 90 & 2825\\
11 & 1269 & 31 & 1619 &51 & 2127 &71 & 2517 & 91 & 2831\\
12 & 1298 & 32 & 1641 & 52 & 2128 & 72 & 2529 & 92 & 2859\\
13 & 1325 & 33 & 1678 & 53 & 2133 & 73 & 2532 & 93 & 2861\\
14 & 1329 & 34 & 1728 & 54 & 2144 & 74 & 2548 & 94 & 2864\\
15 & 1342 & 35 & 1732 & 55 & 2164 & 75 & 2568 & 95 & 2894\\
16 & 1374 & 36 & 1800 & 56 & 2199 & 76 & 2614 & 96 & 2897\\
17 & 1377 & 37 & 1809 & 57 & 2210 & 77 & 2624 & 97 & 2947\\
18 & 1381 & 38 & 1842 & 58 & 2221 & 78 & 2638 & 98 & 2963\\
19 & 1434 & 39 & 1851 & 59 & 2300 & 79 & 2639 & 99 & 2973\\
20 & 1464 & 40 & 1873 & 60 & 2303 & 80 & 2695 & 100 & 2984\\\hline
    \end{tabular}}
    \caption{The individual demands for the 100 involved buyers.}
    \label{tab:demandsexperimental}
\end{table}\normalsize

All buyers are assumed to operate under the same Schotanus price function $P$, see (\ref{eq:psch}). Thus, we keep its parameters $a$, $b$ and $\eta$ fixed across all agents and let the differences in individual demand drive the heterogeneity in the resulting outcomes. The parameters $a$ and $b$ were also randomly generated before constructing the instance according to $a\sim U^d(0,10)$ and $b\sim U^d(0,10)$.\footnote{$U^d(a,b)$ denotes a discrete uniform distribution over the integers in the interval $[a,b]$, with each possible value having equal probability.} Thus, we consider $a=6$ and $b=2$. For the parameter $\eta$, following Assumption 3 in \cite{schotanus2008unfair}, we set $\eta=1/2$ as a reasonable representative value in practice for the elasticity parameter. Since the price depends on the size of the group, buyers' purchasing decisions are interdependent. This provides a simple benchmark to study how differences in individual demand affect group participation and the resulting group configuration.

Once the multi-agent group-buying situation $(N,q,P)$ is specified, we compute the different cost allocations discussed throughout this paper. In this setting, the aggregate quantity demanded by the group is $q_N=205890$, which, according to the pricing function $P$, corresponds to a unit price of $P(q_N)=6.005$. Based on these values, we then determine the associated cost structure, obtaining non-separable costs of $NSC(N,q,P)=452.518$ and total cost $C^{GB}(N,q,P)=1236369$.

Table \ref{tab:resultscostgap} in Appendix \ref{app:numcost-gap} presents the corresponding cost-gap-based allocations for the 100 buyers involved. The numerical results show, first, that the three allocation rules considered, $EP$, $ENSC$, and $SCRB$, reduce the cost incurred by each buyer compared with the cost they would incur when purchasing individually, $C^{GB}({i},q,P)$. Moreover, the allocations obtained under the three rules lie between the separable cost, $SC(N,q,P)$, and the individual cost of each buyer. The differences among the three rules are relatively small throughout the sample. For lower demand levels, $ENSC$ generally yields the highest allocation, whereas $EP$ and $SCRB$ produce very similar values. As demand increases, the differences among the three allocation rules become relatively smaller, and $EP$, $ENSC$, and $SCRB$ converge closely to one another. Overall, the numerical results indicate that, although the three rules distribute the cost generated by group buying somewhat differently, they lead to fairly similar cost allocations and ensure that all buyers benefit from participating in the group purchase compared with purchasing individually.

To complete our analysis, we next consider the bankruptcy-driven approach used throughout the paper, for which the bankruptcy problem specified in Definition \ref{defBANK} is required. For this example, we obtain $E^{GB}=8069.61$, and Table \ref{tab:resultscostbank} in Appendix \ref{app:numBD}  reports the individual costs and the corresponding claims $d_i^{GB}$, together with the resulting cost allocations under the proportional rule ($\phi^{PROP,GB}$), the constrained equal awards rule ($\phi^{CEA,GB}$), the constrained equal losses rule ($\phi^{CEL,GB}$), and the Talmud rule ($\phi^{TAL,GB}$) for the corresponding bankruptcy problems. The results show that all four bankruptcy-driven rules lead to substantial reductions in the costs borne by the buyers relative to their individual costs, $C^{GB}(\{i\},q,P)$. 

\begin{figure}[h!]
\par\medskip
\noindent
\begin{subfigure}[b]{0.3\textwidth}
\begin{tikzpicture}
\begin{axis}[
    width=6.5cm,
    height=5.5cm,
    xlabel={Buyers},
    ylabel={},
    xtick={0,20,40,60,80,100},
    xmin=-3, xmax=104,
    ymin=-1.5, ymax=1.5,
    axis lines=left,
]
\addplot[black, thin, no marks]
    table[x=i,y=diffEPprop] {\allocations};

\addplot[black, dashed, thin, no marks]
    coordinates {(-3,0) (104,0)};
\end{axis}
\end{tikzpicture}
\caption{$EP_i-\phi_i^{PROP,GB}$}
\label{fig:ep-prop}
\end{subfigure}
\hfill
\begin{subfigure}[b]{0.3\textwidth}
\begin{tikzpicture}
\begin{axis}[
    width=6.5cm,
    height=5.5cm,
    xlabel={Buyers},
    ylabel={},
    xtick={0,20,40,60,80,100},
    xmin=-3, xmax=104,
    ymin=-10.5, ymax=5.5,
    axis lines=left,
]
\addplot[black, thin, no marks]
    table[x=i,y=diffEPcea] {\allocations};

\addplot[black, dashed, thin, no marks]
    coordinates {(-3,0) (104,0)};
\end{axis}
\end{tikzpicture}
\caption{$EP_i-\phi_i^{CEA,GB}$}
\label{fig:ep-cea}
\end{subfigure}
\hfill
\begin{subfigure}[b]{0.3\textwidth}
\begin{tikzpicture}
\begin{axis}[
    width=6.5cm,
    height=5.5cm,
    xlabel={Buyers},
    ylabel={},
    xtick={0,20,40,60,80,100},
    xmin=-3, xmax=104,
    ymin=-3, ymax=3,
    axis lines=left,
]
\addplot[black, thin, no marks]
    table[x=i,y=diffEPcel] {\allocations};

\addplot[black, dashed, thin, no marks]
    coordinates {(-3,0) (104,0)};
\end{axis}
\end{tikzpicture}
\caption{$EP_i-\phi_i^{CEL,GB}$}
\label{fig:ep-cel}
\end{subfigure}
\begin{subfigure}[b]{0.3\textwidth}
\begin{tikzpicture}
\begin{axis}[
    width=6.5cm,
    height=5.5cm,
    xlabel={Buyers},
    ylabel={ },
    xtick={0,20,40,60,80,100},
    xmin=-3, xmax=104,
    ymin=-12, ymax=6,
    axis lines=left,
    legend pos=north east,
]
\addplot[black, thin, no marks]
    table[x=i,y=diffPC] {\allocations};

\addplot[black, dashed, thin, no marks]
    coordinates {(-3,0) (104,0)};
\end{axis}
\end{tikzpicture}
\caption{$\phi_i^{PROP,GB}-\phi_i^{CEA,GB}$}
\label{fig:prop-cea}
\end{subfigure}\hfill\hfill
\begin{subfigure}[b]{0.3\textwidth}
\centering
\begin{tikzpicture}
\begin{axis}[
    width=6.5cm,
    height=5.5cm,
    xlabel={Buyers},
    ylabel={ },
    xtick={0,20,40,60,80,100},
    xmin=-3, xmax=104,
    ymin=-2, ymax=2,
    axis lines=left
]
\addplot[black, thin, no marks]
    table[x=i,y=diffPT] {\allocations};

\addplot[black, dashed, thin, no marks]
    coordinates {(-3,0) (104,0)};
\end{axis}
\end{tikzpicture}
\caption{$\phi_i^{PROP,GB}-\phi_i^{CEL,GB}$}
\label{fig:prop-cel}
\end{subfigure}\hfill\hfill
\begin{subfigure}[b]{0.3\textwidth}
\centering
\begin{tikzpicture}
\begin{axis}[
    width=6.5cm,
    height=5.5cm,
    xlabel={Buyers},
    ylabel={ },
    xtick={0,20,40,60,80,100},
    xmin=-3, xmax=104,
    ymin=-6, ymax=13,
    axis lines=left,
    legend pos=north west,
]
\addplot[black, thin, no marks]
    table[x=i,y=diffCT] {\allocations};

\addplot[black, dashed, thin, no marks]
    coordinates {(-3,0) (104,0)};
\end{axis}
\end{tikzpicture}
\caption{$\phi_i^{CEA,GB}-\phi_i^{CEL,GB}$}
\label{fig:cea-cel}
\end{subfigure}

\caption{Pairwise differences between bankruptcy-driven
cost allocations. Positive values indicate that the first
rule assigns a higher cost than the second rule.}
\label{fig:pairwise-differences}

\end{figure}
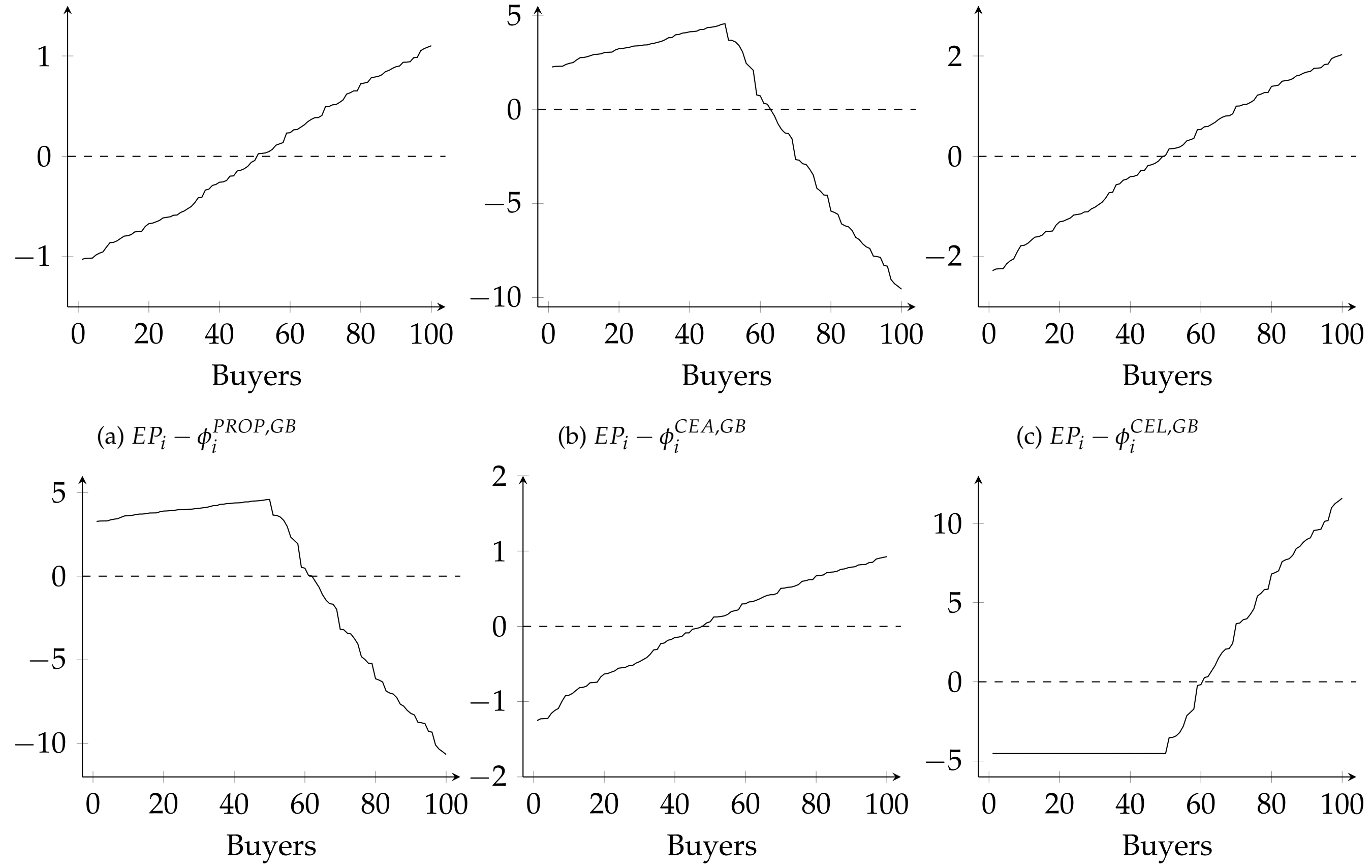

The allocation patterns differ systematically across the rules. For buyers with lower demands in the reported instance, $\phi^{CEA,GB}$ yields the lowest allocation, while $\phi^{CEL,GB}$ and $\phi^{TAL,GB}$ coincide (and also coincide with the ENSC rule) and yield the highest allocation, with $\phi^{PROP,GB}$ lying between them. As demand increases, the relative position of the rules changes. {The differences between the resulting allocations are almost imperceptible in magnitude when plotted on the same scale. For this reason, Figure \ref{fig:pairwise-differences} displays the corresponding pairwise differences between the allocation rules, making their relative differences and how they evolve with demand more clearly visible.}

These results illustrate how the bankruptcy-based interpretation leads to different distributions of the group-buying costs depending on the principle used to allocate the available amount among the buyers. Finally, we can also compare our bankruptcy-driven cost allocations with the cost-gap-based allocations reported in Table \ref{tab:resultscostgap}. In this particular case, we find that $\phi^{PROP,GB}(N,q,P)=SCRB(N,q,P)$, while $\phi^{CEL,GB}=\phi^{TAL,GB}=ENSC(N,q,P)$.

The numerical analysis also highlights the role played by demand heterogeneity in the distribution of cooperative gains. Buyers with different demand levels do not benefit from group purchasing in the same way under the different allocation principles. Hence, heterogeneity in demand has a twofold effect: it determines the magnitude of the savings generated by cooperation and, at the same time, influences the allocation of those savings through the induced claim structure.
The results can also be interpreted through the induced pricing mechanism. Although all buyers face the same underlying quantity discount function, each bankruptcy-driven rule generates a different effective discount factor for each buyer. Therefore, the allocation rule acts as an endogenous pricing mechanism that transforms the savings created by collective purchasing into buyer-specific discounts.
From a managerial perspective, this distinction may be relevant when designing a purchasing group. When buyers differ substantially in their demand levels, the choice of the allocation rule may affect how the cooperative gains are perceived and distributed across participants. Thus, the design of a group-buying mechanism should take into account not only the efficiency gains generated by demand aggregation, but also the distributional consequences of the selected sharing rule.

Another relevant feature of the numerical example is that some allocation rules may coincide for particular demand configurations while diverging for others. The coincidence observed between CEL and Talmud in this instance should therefore be interpreted as a consequence of the underlying claim structure rather than as a general equivalence between the two rules. This is consistent with the analytical results in Section \ref{sec35}, where some relationships among allocation mechanisms are shown to hold generally, while others depend on the particular bankruptcy regime generated by the group-buying situation.

The preceding numerical example also illustrates the practical application of the variable discount price functions introduced above. As discussed, different cost-allocation rules induce different effective prices for individual buyers, reflecting how the gains from collective purchasing are distributed across the group. Expressing these effective prices in the form of the variable discount price function defined in \eqref{generalpricef} allows us to compute the corresponding variable discount factors $b(q_i,q_{-i})$ for each buyer under each allocation rule. These factors provide an alternative perspective on the numerical results, allowing us to examine how the allocation of cooperative gains translates into buyer-specific discounts and how these discounts vary with individual demand and the demand profile of the remaining buyers. Table \ref{tab:b} in Appendix \ref{app:discount} reports the resulting variable discount factors under the Equal Price (EP) rule and the bankruptcy-based allocation rules: Constrained Equal Awards (CEA), Constrained Equal Losses (CEL), Proportional (PROP), and Talmud.

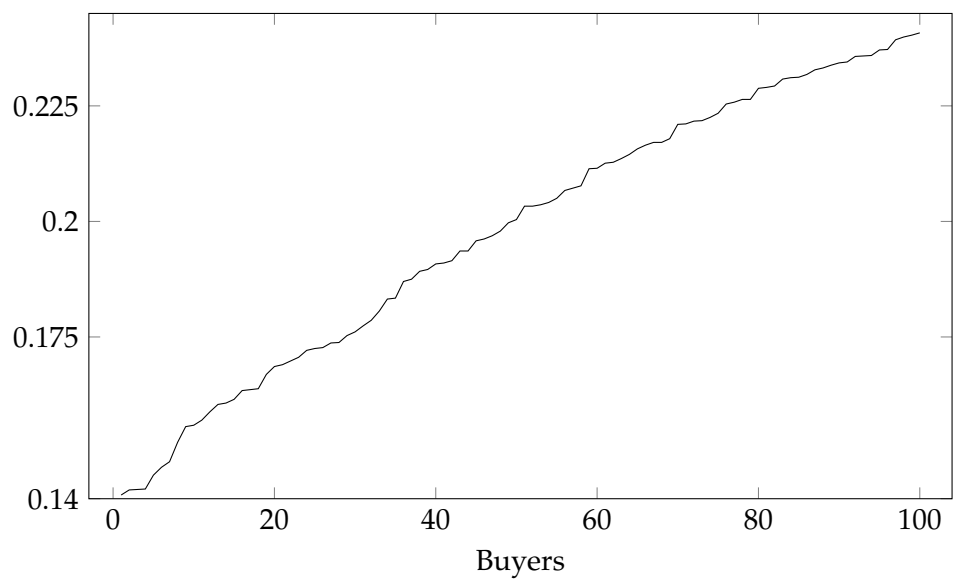
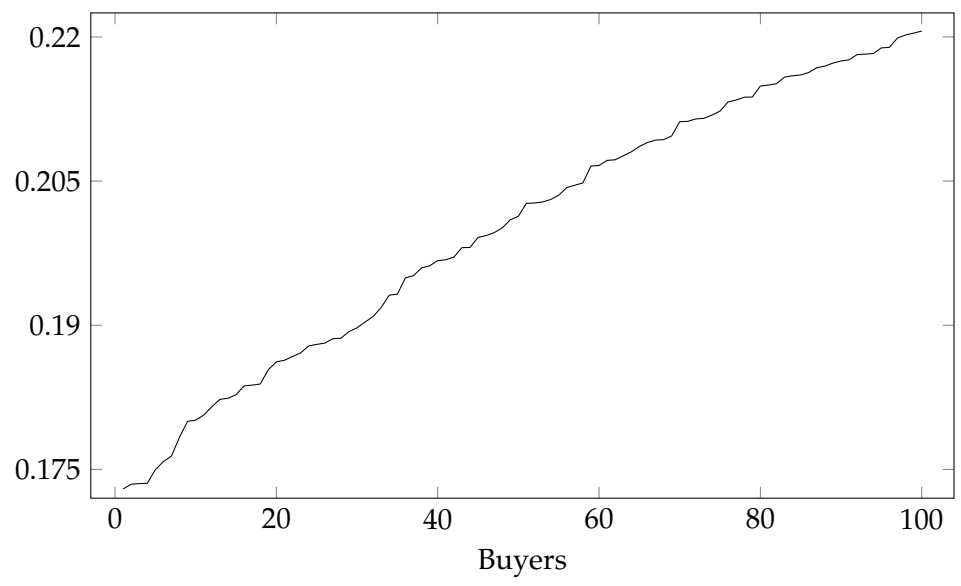
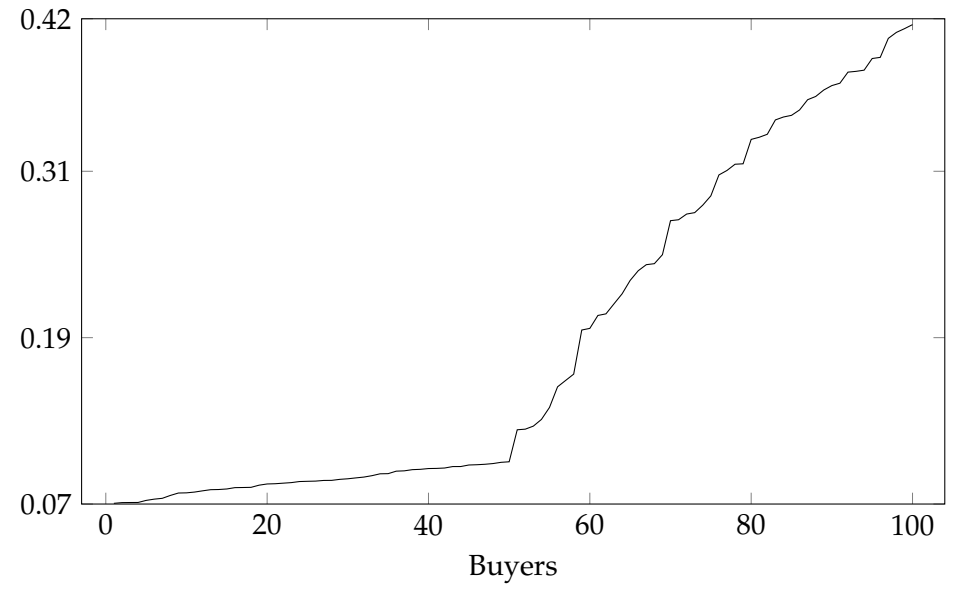
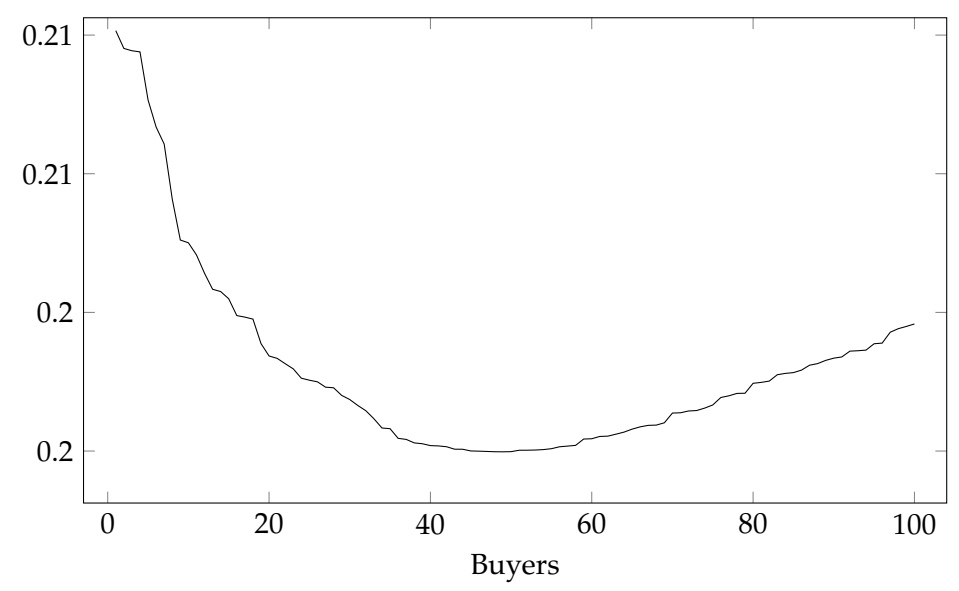
\begin{figure}[!h]
\centering
\begin{subfigure}[b]{0.4\textwidth}
\centering\begin{tikzpicture}[scale=0.5]
\begin{axis}[
    width=13cm,
    height=8cm,
    xlabel={Buyers},
    ylabel={},
    xtick={0,20,40,60,80,100},
    ytick={0.14,0.175,0.200,0.225},
    yticklabel style={
        /pgf/number format/fixed,
        /pgf/number format/precision=3
    },
    xmin=-3,
    xmax=104,
    ymin=0.140,
    ymax=0.245,
]

\addplot[
    black,
    thin,
    no marks
] coordinates {
(1,0.1408)
(2,0.1419)
(3,0.1420)
(4,0.1421)
(5,0.1451)
(6,0.1468)
(7,0.1480)
(8,0.1522)
(9,0.1556)
(10,0.1559)
(11,0.1570)
(12,0.1588)
(13,0.1604)
(14,0.1607)
(15,0.1615)
(16,0.1634)
(17,0.1636)
(18,0.1638)
(19,0.1669)
(20,0.1686)
(21,0.1690)
(22,0.1698)
(23,0.1706)
(24,0.1721)
(25,0.1725)
(26,0.1727)
(27,0.1737)
(28,0.1738)
(29,0.1753)
(30,0.1761)
(31,0.1774)
(32,0.1786)
(33,0.1806)
(34,0.1832)
(35,0.1834)
(36,0.1870)
(37,0.1875)
(38,0.1892)
(39,0.1896)
(40,0.1908)
(41,0.1910)
(42,0.1915)
(43,0.1936)
(44,0.1936)
(45,0.1958)
(46,0.1962)
(47,0.1969)
(48,0.1979)
(49,0.1997)
(50,0.2004)
(51,0.2033)
(52,0.2033)
(53,0.2036)
(54,0.2041)
(55,0.2050)
(56,0.2067)
(57,0.2072)
(58,0.2077)
(59,0.2114)
(60,0.2115)
(61,0.2126)
(62,0.2128)
(63,0.2136)
(64,0.2145)
(65,0.2157)
(66,0.2165)
(67,0.2171)
(68,0.2171)
(69,0.2179)
(70,0.2210)
(71,0.2211)
(72,0.2217)
(73,0.2218)
(74,0.2225)
(75,0.2234)
(76,0.2254)
(77,0.2258)
(78,0.2264)
(79,0.2264)
(80,0.2288)
(81,0.2290)
(82,0.2293)
(83,0.2308)
(84,0.2311)
(85,0.2312)
(86,0.2318)
(87,0.2328)
(88,0.2332)
(89,0.2338)
(90,0.2343)
(91,0.2345)
(92,0.2357)
(93,0.2358)
(94,0.2359)
(95,0.2371)
(96,0.2372)
(97,0.2393)
(98,0.2399)
(99,0.2403)
(100,0.2408)
};

\end{axis}

\end{tikzpicture}
\caption[Ex1]%
{{\small $b^{EP}(q_i,q_{-i})$.}}    
\end{subfigure}\hfil\begin{subfigure}[b]{0.4\textwidth}
\centering
\begin{tikzpicture}[scale=0.5]
\begin{axis}[
    width=13cm,
    height=8cm,
    xlabel={Buyers},
    ylabel={},
    xtick={0,20,40,60,80,100},
    ytick={0.175,0.190,0.205,0.22},
    yticklabel style={
        /pgf/number format/fixed,
        /pgf/number format/precision=3
    },
    xmin=-3,
    xmax=104,
    ymin=0.172,
    ymax=0.2225,
]

\addplot[
    black,
    thin,
    no marks
] coordinates {
(1, 0.17296183)
(2, 0.17345170)
(3, 0.17351675)
(4, 0.17354925)
(5, 0.17496425)
(6, 0.17581846)
(7, 0.17638222)
(8, 0.17835141)
(9, 0.18000220)
(10, 0.18012084)
(11, 0.18065239)
(12, 0.18150098)
(13, 0.18228265)
(14, 0.18239778)
(15, 0.18277077)
(16, 0.18368135)
(17, 0.18376618)
(18, 0.18387914)
(19, 0.18536084)
(20, 0.18618755)
(21, 0.18635188)
(22, 0.18673403)
(23, 0.18711441)
(24, 0.18784306)
(25, 0.18800412)
(26, 0.18813809)
(27, 0.18859200)
(28, 0.18864524)
(29, 0.18933428)
(30, 0.18972927)
(31, 0.19033136)
(32, 0.19090333)
(33, 0.19185682)
(34, 0.19312896)
(35, 0.19322994)
(36, 0.19492934)
(37, 0.19515187)
(38, 0.19596317)
(39, 0.19618319)
(40, 0.19671880)
(41, 0.19681585)
(42, 0.19708221)
(43, 0.19806834)
(44, 0.19809226)
(45, 0.19913886)
(46, 0.19932791)
(47, 0.19963431)
(48, 0.20012721)
(49, 0.20096632)
(50, 0.20133693)
(51, 0.20269150)
(52, 0.20271430)
(53, 0.20282822)
(54, 0.20307836)
(55, 0.20353156)
(56, 0.20431973)
(57, 0.20456616)
(58, 0.20481199)
(59, 0.20656020)
(60, 0.20662600)
(61, 0.20715090)
(62, 0.20721632)
(63, 0.20762972)
(64, 0.20804148)
(65, 0.20860232)
(66, 0.20901028)
(67, 0.20926712)
(68, 0.20930987)
(69, 0.20969381)
(70, 0.21117389)
(71, 0.21121588)
(72, 0.21146747)
(73, 0.21153027)
(74, 0.21186463)
(75, 0.21228112)
(76, 0.21323305)
(77, 0.21343891)
(78, 0.21372646)
(79, 0.21374697)
(80, 0.21488954)
(81, 0.21499099)
(82, 0.21513286)
(83, 0.21581940)
(84, 0.21596023)
(85, 0.21604063)
(86, 0.21630152)
(87, 0.21680154)
(88, 0.21696109)
(89, 0.21727951)
(90, 0.21749791)
(91, 0.21761686)
(92, 0.21817034)
(93, 0.21820978)
(94, 0.21826890)
(95, 0.21885848)
(96, 0.21891727)
(97, 0.21989281)
(98, 0.22020327)
(99, 0.22039689)
(100, 0.22060950)
};

\end{axis}

\end{tikzpicture}
\caption[Ex1]%
{{\small $b^{PROP}(q_i,q_{-i})$.}}    
\end{subfigure}

\begin{subfigure}[b]{0.4\textwidth}
\centering
\begin{tikzpicture}[scale=0.5]
\begin{axis}[
    width=13cm,
    height=8cm,
    xlabel={Buyers},
    ylabel={},
    xtick={0,20,40,60,80,100},
    ytick={0.07,0.190,0.310,0.420},
    yticklabel style={
        /pgf/number format/fixed,
        /pgf/number format/precision=3
    },
    xmin=-3,
    xmax=104,
    ymin=0.07,
    ymax=0.42,
]

\addplot[
    black,
    thin,
    no marks
] coordinates {
(1, 0.07050739)
(2, 0.07102473)
(3, 0.07109343)
(4, 0.07112775)
(5, 0.07262211)
(6, 0.07352421)
(7, 0.07411959)
(8, 0.07619921)
(9, 0.07794256)
(10, 0.07806785)
(11, 0.07862921)
(12, 0.07952539)
(13, 0.08035090)
(14, 0.08047248)
(15, 0.08086639)
(16, 0.08182803)
(17, 0.08191762)
(18, 0.08203691)
(19, 0.08360170)
(20, 0.08447477)
(21, 0.08464831)
(22, 0.08505190)
(23, 0.08545361)
(24, 0.08622312)
(25, 0.08639320)
(26, 0.08653469)
(27, 0.08701405)
(28, 0.08707027)
(29, 0.08779796)
(30, 0.08821509)
(31, 0.08885095)
(32, 0.08945500)
(33, 0.09046195)
(34, 0.09180543)
(35, 0.09191207)
(36, 0.09370677)
(37, 0.09394178)
(38, 0.09479858)
(39, 0.09503093)
(40, 0.09559658)
(41, 0.09569907)
(42, 0.09598037)
(43, 0.09702180)
(44, 0.09704706)
(45, 0.09815234)
(46, 0.09835200)
(47, 0.09867558)
(48, 0.09919612)
(49, 0.10008229)
(50, 0.10047368)
(51, 0.12356156)
(52, 0.12400250)
(53, 0.12620257)
(54, 0.13101560)
(55, 0.13967236)
(56, 0.15453653)
(57, 0.15913504)
(58, 0.16369935)
(59, 0.19551137)
(60, 0.19668706)
(61, 0.20601059)
(62, 0.20716589)
(63, 0.21443132)
(64, 0.22160914)
(65, 0.23129259)
(66, 0.23826977)
(67, 0.24263410)
(68, 0.24335834)
(69, 0.24983645)
(70, 0.27436448)
(71, 0.27505020)
(72, 0.27914748)
(73, 0.28016724)
(74, 0.28557554)
(75, 0.29226469)
(76, 0.30735737)
(77, 0.31058575)
(78, 0.31507462)
(79, 0.31539389)
(80, 0.33298817)
(81, 0.33453241)
(82, 0.33668716)
(83, 0.34703541)
(84, 0.34914206)
(85, 0.35034225)
(86, 0.35422485)
(87, 0.36161499)
(88, 0.36395889)
(89, 0.36861667)
(90, 0.37179591)
(91, 0.37352223)
(92, 0.38150639)
(93, 0.38207220)
(94, 0.38291980)
(95, 0.39132318)
(96, 0.39215633)
(97, 0.40585435)
(98, 0.41016431)
(99, 0.41284036)
(100, 0.41576845)
};

\end{axis}

\end{tikzpicture}
\caption[Ex2]%
{{\small  $b^{CEA}(q_i,q_{-i})$.}}    
\end{subfigure}\hfil\begin{subfigure}[b]{0.4\textwidth}
\centering
\begin{tikzpicture}[scale=0.5]
\begin{axis}[
    width=13cm,
    height=8cm,
    xlabel={Buyers},
    ylabel={ },
    xtick={0,20,40,60,80,100},
    ytick={0.200,0.204,0.208,0.212},
    yticklabel style={
        /pgf/number format/fixed,
        /pgf/number format/precision=2
    },
    xmin=-3,
    xmax=104,
    ymin=0.1985,
    ymax=0.2125,
]

\addplot[
    black,
    thin,
    no marks
] coordinates {
(1, 0.21212682)
(2, 0.21161519)
(3, 0.21154838)
(4, 0.21151509)
(5, 0.21012805)
(6, 0.20934750)
(7, 0.20885482)
(8, 0.20726748)
(9, 0.20608781)
(10, 0.20600806)
(11, 0.20565881)
(12, 0.20512792)
(13, 0.20466711)
(14, 0.20460147)
(15, 0.20439270)
(16, 0.20390742)
(17, 0.20386395)
(18, 0.20380651)
(19, 0.20309984)
(20, 0.20274220)
(21, 0.20267414)
(22, 0.20251968)
(23, 0.20237118)
(24, 0.20210107)
(25, 0.20204387)
(26, 0.20199697)
(27, 0.20184263)
(28, 0.20182498)
(29, 0.20160517)
(30, 0.20148623)
(31, 0.20131463)
(32, 0.20116227)
(33, 0.20093077)
(34, 0.20066430)
(35, 0.20064517)
(36, 0.20036625)
(37, 0.20033560)
(38, 0.20023506)
(39, 0.20021077)
(40, 0.20015688)
(41, 0.20014790)
(42, 0.20012448)
(43, 0.20005320)
(44, 0.20005177)
(45, 0.20000260)
(46, 0.19999650)
(47, 0.19998837)
(48, 0.19997985)
(49, 0.19997804)
(50, 0.19998224)
(51, 0.20002295)
(52, 0.20002397)
(53, 0.20002923)
(54, 0.20004173)
(55, 0.20006769)
(56, 0.20012281)
(57, 0.20014260)
(58, 0.20016355)
(59, 0.20034635)
(60, 0.20035436)
(61, 0.20042117)
(62, 0.20042985)
(63, 0.20048652)
(64, 0.20054605)
(65, 0.20063199)
(66, 0.20069799)
(67, 0.20074102)
(68, 0.20074830)
(69, 0.20081503)
(70, 0.20109554)
(71, 0.20110402)
(72, 0.20115548)
(73, 0.20116848)
(74, 0.20123878)
(75, 0.20132886)
(76, 0.20154498)
(77, 0.20159356)
(78, 0.20166252)
(79, 0.20166748)
(80, 0.20195412)
(81, 0.20198051)
(82, 0.20201767)
(83, 0.20220160)
(84, 0.20224017)
(85, 0.20226232)
(86, 0.20233482)
(87, 0.20247645)
(88, 0.20252238)
(89, 0.20261510)
(90, 0.20267951)
(91, 0.20271486)
(92, 0.20288189)
(93, 0.20289395)
(94, 0.20291207)
(95, 0.20309530)
(96, 0.20311382)
(97, 0.20342780)
(98, 0.20353030)
(99, 0.20359485)
(100, 0.20366627)
};

\end{axis}

\end{tikzpicture}
\caption[Ex2]%
{{\small $b^{CEL}(q_i,q_{-i})$ and $b^{TAL}(q_i,q_{-i})$.}}    
\end{subfigure}

\caption[ Ex ]
{\small Evolution of the variable discount factors by buyer.} \label{im:b}
\end{figure}

\textcolor{black}{Figure \ref{im:b}  complements the analysis by illustrating how the variable discount factors evolve across buyers under the different allocation rules. Under the EP and PROP rules, these factors increase rapidly from the first buyer onwards. In contrast, under the CEA rule, they increase much more slowly for the initial buyers, with the rate of growth becoming more pronounced around the 50th buyer. The CEL and Talmud rules, that coincide in this example, exhibit a different pattern. In both cases, the variable discount factors initially decrease, reaching their minimum values around  the 50th buyer, before increasing again towards the final buyers.}

\bigskip

{\color{black}The numerical example above illustrates how the choice of cost-allocation rule affects the behavior of variable discount factors across buyers. In order to assess the robustness of these findings, we conduct an extensive simulation study covering a broad range of group-buying scenarios. Specifically, we consider different group sizes, values of $\eta$, and demand distributions. Appendix \ref{appD} reports the simulation results for the demand distribution considered in this section. We also consider alternative demand distributions, including the log-normal and normal distributions, and obtain analogous conclusions. Across all these scenarios, the qualitative patterns observed in the reference numerical example persist, suggesting that our findings are robust to changes in group size, parameter values, and demand assumptions. For the sake of brevity, we do not report the simulation results for the normal and log-normal distributions.}

\section{Conclusions and Future Research}\label{sec5}

This paper develops an integrated framework for allocating costs and cooperative savings in group-buying environments with quantity discounts. Although demand aggregation generates economies of scale, the resulting cost savings do not determine how these gains should be distributed among participating buyers. The allocation problem therefore involves a choice among alternative distributional principles, with implications for both individual payments and the effective prices faced by buyers with heterogeneous demands.

We address this problem through two complementary approaches. The first relies on cost-gap allocation methods, which decompose the total group-buying cost into separable and non-separable components. The Equal Price (EP), Egalitarian Non-Separable Cost (ENSC), and Surplus Cost Ratio Based (SCRB) rules provide alternative ways of distributing the costs associated with cooperation. These mechanisms preserve the cost advantage of collective purchasing while producing different distributions of individual payments, illustrating how the choice of allocation principle affects the division of the gains from demand aggregation.

The second, and main, contribution is to formulate group-buying cost allocation as a bankruptcy problem \citep{o1982problem}. The aggregate savings generated by cooperation define the estate, while the differences between buyers' stand-alone costs and their separable costs determine their respective claims. Bankruptcy theory thus provides an endogenous allocation framework grounded in the economic structure of group purchasing rather than an external set of fairness criteria imposed on the problem. This formulation enables established bankruptcy rules to be applied to the allocation of cooperative savings and provides a common basis for comparing alternative mechanisms.

We examine the Proportional (PROP), Constrained Equal Awards (CEA), Constrained Equal Losses (CEL), and Talmud rules. Because these rules embody different principles for distributing a limited estate among competing claims, they may generate substantially different individual payments even when applied to the same purchasing group. The analysis also identifies structural relationships between bankruptcy-driven and cost-gap-based mechanisms. In particular, the proportional bankruptcy allocation coincides with the SCRB allocation, while the CEL-based allocation coincides with ENSC under specific conditions. Beyond these cases, the two approaches may produce different distributions of cooperative gains, with no general ordering of individual buyers' allocations. These relationships clarify both the connections between established allocation methods and the circumstances in which their distributional implications diverge. A further contribution is to establish a pricing interpretation of the bankruptcy-based allocations. Building on the quantity-discount specification of \cite{schotanus2009unraveling}, we show how each allocation rule induces an individual price function of the same general form, with a variable discount factor that depends on the buyer's own demand and the demand profile of the remaining group members. This representation links the ex post allocation of cooperative savings to the effective prices paid by individual buyers. It also provides a way to examine how alternative distributional principles translate into buyer-specific discounts, thereby connecting cost-sharing design with procurement pricing. {As an alternative to the approach pursued in this paper, one could also consider bankrupcty-driven allocation rules of the form
\begin{equation*}
\psi_i^{\chi,GB}(N,q,P)=SC_i(N,q,P)+\chi_i(NSC(N,q,P),d^{GB}), \mbox{ for every }i \in N,
\end{equation*}
where $\chi$ is a bankruptcy rule. Since $(E^{GB},d^{GB})$ and $(NSC(N,q,P),d^{GB})$ define dual bankruptcy problems \citep{thomson2015axiomatic}, such an approach would be closely related to the allocation rules proposed here, although the precise correspondence would depend on the bankruptcy rule considered.}

{The numerical example and simulation study complement the analytical results by illustrating how allocation rules shape the distribution of costs and the behavior of the associated variable discount factors. The simulation results show that cost allocations exhibit increasing behavior across buyers under all allocation rules, regardless of group size or the parameter regime considered. In contrast, the behavior of variable discount factors depends on the allocation rule: EP, PROP, and CEA display increasing patterns, whereas CEL and TAL exhibit non-monotonic behavior in all simulated scenarios. These patterns remain unchanged across the different number of buyers and parameter regimes for $\eta$ considered. Moreover, analogous conclusions are obtained when alternative demand distributions, such as the normal and log-normal distributions, are used. Taken together, the numerical and simulation analyses reinforce the conclusion that the choice of allocation rule plays a crucial role in determining how effective discounts are distributed across buyers, even when the underlying cost allocations exhibit the same qualitative behavior.}

Several directions for future research follow from this framework. Other bankruptcy rules could be examined to determine how their axiomatic properties (including monotonicity, consistency, solidarity, and invariance under changes in the demand profile) translate into properties of group-buying allocations and induced price functions. A further extension would endogenize group formation. In the present analysis, the purchasing group is given, whereas in practice buyers may decide whether to join or leave a collective purchasing arrangement. Allowing participation to depend on individual payments and effective prices would make it possible to study coalition formation, participation incentives, and the stability of group-buying arrangements under alternative allocation rules. This would establish a direct link between bankruptcy-based cost allocation and cooperative game-theoretic approaches to stability. Finally, the pricing framework could be extended to alternative discount schedules, multiple suppliers, capacity constraints, heterogeneous products, and uncertain demand. Such extensions would help establish the robustness and practical applicability of the proposed framework in more complex procurement settings.

Overall, this paper connects group buying, bankruptcy theory, and endogenous pricing within a unified analytical framework. By linking the economic gains generated by demand aggregation to alternative allocation principles and their corresponding buyer-specific prices, it provides a basis for evaluating the distributional consequences of collective purchasing. The framework can support the comparison and design of cost-allocation mechanisms in settings where efficiency gains must be reconciled with heterogeneous buyer interests and competing notions of fairness.

\subsection*{Acknowledgments}
This work is part of the R+D+I project grants PID2025-170122NB-C33, PID2022-137211NB-I00 and PID2025-170122NB-C31, that were funded by MCIN/AEI/10.13039/501100011033/ and by
"ERDF A way of making Europe"/EU. This research was also partially supported by the CIPROM/2024/34 grant, funded by the Conselleria de Educaci\'on, Cultura, Universidades y Empleo, Generalitat Valenciana, and by the ED431C 2025/03  and ED431C 2024/26, grants, funded by \emph{Conseller\'ia de Educaci\'on, Ciencia, Universidades e Formaci\'on Profesional, Xunta de Galicia}

\subsection*{Data and Code availability}
The R code of all the allocation rules presented in this paper, as well as the used code in the simulation study, are available at \url{https://github.com/AlejandroSaavedraNieves/Group-Buying-.git}.

\bibliographystyle{apalike}
\bibliography{referencias}

\clearpage

\appendix

\setcounter{table}{0}
\renewcommand{\thetable}{\Alph{section}.\arabic{table}}
\section{Additional results}\label{app0}

This section presents additional results concerning the non-separable costs and the associated bankruptcy problem.

\subsection{The non-separable costs}\label{appNSC}

Take $(N,q,P)\in GB^N$, where $P$ is a Schotanus price function, and consider the non-separable costs $NSC(N,q,P)$ defined in Equation~(\ref{NSCformula}). Recall that these costs are given by
\begin{equation*}
    \begin{aligned}
 NSC(N,q,P)
 =& b\bigg(q_N^{1-\eta}-\sum_{i\in N}(q_N^{1-\eta}-q_{N\setminus \{i\}}^{1-\eta}  ) \bigg).
 \end{aligned}
\end{equation*}
We now prove that these costs are non-negative. Since $q_N>0$ and
$q_{N\setminus\{i\}}=q_N-q_i
=q_N\left(1-\frac{q_i}{q_N}\right)$, we have
\[
q_{N\setminus\{i\}}^{1-\eta}
=q_N^{1-\eta}
\left(1-\frac{q_i}{q_N}\right)^{1-\eta}.
\]
Consequently,
\[
\sum_{i\in N}q_{N\setminus\{i\}}^{1-\eta}
=q_N^{1-\eta}
\sum_{i\in N}\left(1-\frac{q_i}{q_N}\right)^{1-\eta}.
\]

Since $0<1-\eta<1$, we have $x^{1-\eta}\geq x$ for all $x\in[0,1]$. Therefore,\[
\left(1-\frac{q_i}{q_N}\right)^{1-\eta}
\ge
1-\frac{q_i}{q_N}.
\]
Summing over $i\in N$ gives
\[
\sum_{i\in N}\left(1-\frac{q_i}{q_N}\right)^{1-\eta}
\ge
\sum_{i\in N}\left(1-\frac{q_i}{q_N}\right)
=
n-1.
\]

\noindent Multiplying by $q_N^{1-\eta}$ yields $\sum_{i\in N} q_{N\setminus\{i\}}^{1-\eta}
\ge
(n-1)q_N^{1-\eta}$. Hence, since $b\geq 0$, $NSC(N,q,P)\geq 0$.

\subsection{The bankruptcy problem}\label{appBANK}

Take $(N,q,P)\in GB^N$ and consider the associated bankruptcy problem $(E^{GB},d^{GB})$ defined in Definition~\ref{defBANK}. The claims and the estate are given by

\[
d^{GB}_i = b\bigg(q_i^{1-\eta} + q_{N\setminus\{i\}}^{1-\eta} - q_N^{1-\eta}\bigg), \text{ for all } i \in N, \mbox{ and }
E^{GB} = b\bigg(\sum_{i\in N} q_i^{1-\eta} - q_N^{1-\eta}\bigg).
\]
We first prove that $d_i^{GB}\geq 0$ for all $i\in N$. To this end, observe that, for every $\eta\in(0,1)$, the function $f(x)=x^{1-\eta}$ is increasing and concave on $\mathbb{R}_+$. In particular, it is subadditive, that is,
\[
(x+y)^{1-\eta} \le x^{1-\eta} + y^{1-\eta}, \qquad \forall x,y \ge 0.
\]

Since $q_N = q_i + q_{N\setminus\{i\}}$, we obtain $q_N^{1-\eta}
\le
q_i^{1-\eta} + q_{N\setminus\{i\}}^{1-\eta}$, which implies that
\[
d_i^{GB}
=b\bigg(
q_i^{1-\eta} + q_{N\setminus\{i\}}^{1-\eta} - q_N^{1-\eta}\bigg)
\ge 0,
\qquad \forall i \in N.
\]

Next, we show that the estate $E^{GB}$ is non-negative. By the subadditivity of $f$, applied repeatedly, we obtain\[
q_N^{1-\eta}
=\left(\sum_{i\in N}q_i\right)^{1-\eta}
\leq \sum_{i\in N}q_i^{1-\eta}.\]
Therefore, since $b\geq 0$, we have $E^{GB}\geq 0$.

Finally, we verify that $(E^{GB},d^{GB})$ is a well-defined bankruptcy problem by showing that $\sum_{i\in N}d_i^{GB}\geq E^{GB}$. Summing over all buyers, we obtain
\[
\sum_{i\in N} d_i^{GB}
=b\bigg(
\sum_{i\in N} q_i^{1-\eta}
+
\sum_{i\in N} q_{N\setminus\{i\}}^{1-\eta}
-
n q_N^{1-\eta}\bigg), 
\]and, thus,
\[
\sum_{i\in N} d_i^{GB} - E^{GB}
=b\bigg(
\sum_{i\in N} q_{N\setminus\{i\}}^{1-\eta}
-
(n-1)q_N^{1-\eta}\bigg),
\]
which coincides with the non-separable costs, $NSC(N,q,P)$. As shown in Appendix (\ref{appNSC}), it holds that $NSC(N,q,P)\geq 0$. Consequently, $\sum_{i\in N}d_i^{GB}\geq E^{GB}$. Hence, all the conditions required for $(E^{GB},d^{GB})$ to constitute a well-defined bankruptcy problem are satisfied.

\section{The CEL and the Talmud rules}\label{app1}
This appendix provides additional details complementing Section \ref{sec3}. More specifically, it describes the obtaining of the CEL and Talmud rules in group-buying using bankruptcy.

Let $(N,q,P)\in GB^N$ be a group-buying multi-agent situation and the associated bankruptcy problem $(E^{GB},d^{GB})\in B^N$ specified in Definition \ref{defBANK}. We assume  $d_1^{GB} \leq d_2^{GB} \leq \cdots \leq d_n^{GB}$.

\subsection{The CEL rule}\label{CELbank}

First, we consider the case of the CEL rule. We derive the explicit expression of $CEL^{GB}(N,q,P)$, distinguishing two cases.  First, we assess $\sum_{i\in N}d^{GB}_i-E^{GB}$ to obtain:
\begin{equation*}
	\begin{array}{lcl}
	\sum_{i\in N}d^{GB}_i-E^{GB}&=&b (\sum_{i\in N} (q_i^{1-\eta}+q_{N\setminus\{i\}}^{1-\eta}-q_N^{1-\eta})-\sum_{i\in N}q_i^{1-\eta}+q_N^{1-\eta})\\
	&=& b(q_N^{1-\eta}+\sum_{i\in N} (q_{N\setminus\{i\}}^{1-\eta}-q_N^{1-\eta}))\geq 0.
	\end{array}
\end{equation*}
Using Equation~(\ref{CELorig2}) we determine $CEA(\sum_{i\in N}d^{GB}_i-E^{GB},d^{GB})$. Thus, we distinguish two scenarios:
\begin{itemize}
\item[(a)] If it holds $d_1^{GB} \geq \frac{\sum_{j\in N}d_j^{GB}-E^{GB}}{n}$, then 
\[
CEA_i(	\sum_{i\in N}d^{GB}_i-E^{GB},d^{GB})=\frac{b}{n}\left(q_N^{1-\eta}+\sum_{j\in N}\left(q_{N\setminus\{j\}}^{1-\eta}-q_N^{1-\eta}\right)\right)
\]
and
\begin{equation}\label{CELaApp}
    CEL_i(	\sum_{i\in N}d^{GB}_i-E^{GB},d^{GB})=b[q_i^{1-\eta} + q_{N\setminus\{i\}}^{1-\eta} - q_N^{1-\eta} -\frac{1}{n}(q_N^{1-\eta}+\sum_{j\in N}(q_{N\setminus\{j\}}^{1-\eta}-q_N^{1-\eta}))]
\end{equation}



\item[(b)] Otherwise, let $k\in \{1,\ldots,n-1\}$ be such that
\begin{equation*}
d_k^{GB} \leq \frac{b}{n-k} \bigg(q_N^{1-\eta}+\sum_{j\in N} (q_{N\setminus\{j\}}^{1-\eta}-q_N^{1-\eta})-\sum_{j=1}^k (q_j^{1-\eta}+q_{N\setminus\{j\}}^{1-\eta}-q_N^{1-\eta})\bigg)< d_{k+1}^{GB}.
\end{equation*}Then, it holds that
\begin{equation*}
	CEA_i(\sum_{i\in N}d^{GB}_i-E^{GB},d^{GB})=
	\begin{cases}
		b(q_i^{1-\eta}+q_{N\setminus\{i\}}^{1-\eta}-q_N^{1-\eta}),&\text{ if }i\leq k.\\
		\frac{b}{n-k}( q_N^{1-\eta}- \sum_{j=k+1}^n (q_N^{1-\eta}-q_{N\setminus\{j\}}^{1-\eta})-\sum_{j=1}^k q_j^{1-\eta}),&\text{ if }k<i\leq n.
	\end{cases}
\end{equation*} 
and, hence, it holds that
\begin{equation}\label{CELbApp}
	CEL_i(E^{GB},d^{GB})=
	\begin{cases}
		0,&\text{ if }i\leq k.\\
	b\bigg(q_i^{1-\eta}+q_{N\setminus\{i\}}^{1-\eta}-q_N^{1-\eta}-	\frac{q_N^{1-\eta}- \sum_{j=k+1}^n (q_N^{1-\eta}-q_{N\setminus\{j\}}^{1-\eta})-\sum_{j=1}^k q_j^{1-\eta}}{n-k}\bigg),&\text{ if  }k<i\leq n.
	\end{cases}
\end{equation}\end{itemize}




\subsection{The Talmud rule}\label{TalmudBank}

This appendix focuses on obtaining the explicit expression of $TAL^{GB}(N,q,P)$. Again, we distinguish two cases. 

\begin{itemize}

\item[(a)] If $E^{GB} \leq \frac{1}{2}\sum_{i\in N} d^{GB}_i$, take $k\in\{0,\ldots,n-1\}$ such that $d_0^{GB}=0$ and 
\begin{equation*}
\frac{d^{GB}_k}{2}
\leq \frac{b}{n-k}\left(\bigg(\sum_{i\in N} q_i^{1-\eta} - q_N^{1-\eta}\bigg)-\sum_{j=1}^k \frac{\bigg(q_j^{1-\eta} + q_{N\setminus\{j\}}^{1-\eta} - q_N^{1-\eta}\bigg)}{2}\right)
<\frac{d^{GB}_{k+1}}{2}.
\end{equation*}

Then,
\begin{equation}\label{TalaApp}
TAL_i(E^{GB},d^{GB})=
\begin{cases}
\frac{b}{2}\left(q_i^{1-\eta}+q_{N\setminus\{i\}}^{1-\eta}-q_N^{1-\eta}\right), & \text{if } i\leq k,\\[6pt]
\frac{b}{n-k}\left(\big(\sum_{i\in N} q_i^{1-\eta} - q_N^{1-\eta})-\sum_{j=1}^k \frac{(q_j^{1-\eta} + q_{N\setminus\{j\}}^{1-\eta} - q_N^{1-\eta})}{2}\right), & \text{if } i>k.
\end{cases}
\end{equation}

\item[(b)] Otherwise, if $E^{GB} \geq \frac{1}{2}\sum_{i\in N} d^{GB}_i$, take $k\in\{0,\ldots,n-1\}$ with $d_0^{GB}=0$ such that
\begin{equation*}
\frac{d^{GB}_k}{2}
\leq \frac{b}{n-k}\left(\sum_{j\in N}\bigg( q_{N\setminus\{j\}}^{1-\eta} - q_N^{1-\eta}\bigg) + q_N^{1-\eta} - \sum_{j=1}^k \frac{\bigg(q_j^{1-\eta} + q_{N\setminus\{j\}}^{1-\eta} - q_N^{1-\eta}\bigg)}{2}\right)
<\frac{d^{GB}_{k+1}}{2}.
\end{equation*}

Then,

\begin{equation}\label{TALbApp}\footnotesize
TAL_i(E^{GB},d^{GB})=
\begin{cases}
\displaystyle
\frac{b\bigg(q_i^{1-\eta}+q_{N\setminus\{i\}}^{1-\eta}-q_N^{1-\eta}\bigg)}{2},
& \text{if } i\leq k,\\[8pt]
\begin{aligned}
&\displaystyle {b\bigg(q_i^{1-\eta}+q_{N\setminus\{i\}}^{1-\eta}-q_N^{1-\eta}\bigg)}-\\[-2pt]
&\displaystyle\quad-\frac{b}{n-k}\left(
\sum_{j\in N}\bigg(q_{N\setminus\{j\}}^{1-\eta}-q_N^{1-\eta}\bigg)
+q_N^{1-\eta}
-\sum_{j=1}^{k}\frac{\bigg(q_j^{1-\eta}+q_{N\setminus\{j\}}^{1-\eta}-q_N^{1-\eta}\bigg)}{2}
\right),
\end{aligned}
& \text{if } i>k.
\end{cases}
\end{equation}

\end{itemize}

\section{Numerical results}\label{app2}

This appendix reports the numerical results for the group-buying situation presented in Section \ref{sec4}. For ease of reference and to improve the readability of the main text, the detailed numerical results are provided here.

\subsection{Cost-gap cost allocations}\label{app:numcost-gap}

Table \ref{tab:resultscostgap} presents the cost-gap-based allocations for the 100 buyers involved in the group-buying situation. For each buyer, the table reports demand, individual cost, separable cost, and the allocations obtained under the three cost-gap-based rules: $EP$, $ENSC$, and $SCRB$.

   \renewcommand{\arraystretch}{0.9}

\footnotesize
\setlength{\LTleft}{\fill}
\setlength{\LTright}{\fill}
\begin{longtable}{c|ccc|ccc}
\caption{Demands, individual costs, separable costs and the cost-gap cost allocations}
\label{tab:resultscostgap}\\
\hline
$i$ & $q_i$ & $C^{GB}(\{i\},q,P)$ & $SC_i$ &
$EP_i$ & $ENSC_i$ & $SCRB_i$\\
\hline
\endfirsthead

\caption[]{Demands, individual costs, separable costs and the cost-gap cost allocations}\\
\hline
$i$ & $q_i$ & $C^{GB}(\{i\},q,P)$ & $SC_i$ &
$EP_i$ & $ENSC_i$ & $SCRB_i$\\
\hline
\endhead
1	&	1021	&	6189.906	&	6128.253	&	6130.5	&	6132.778	&	6131.527	\\
2	&	1036	&	6280.374	&	6218.286	&	6220.566	&	6222.811	&	6221.583	\\
3	&	1038	&	6292.436	&	6230.29	&	6232.575	&	6234.816	&	6233.59	\\
4	&	1039	&	6298.467	&	6236.293	&	6238.58	&	6240.818	&	6239.594	\\
5	&	1083	&	6563.818	&	6500.39	&	6502.774	&	6504.915	&	6503.758	\\
6	&	1110	&	6726.633	&	6662.45	&	6664.893	&	6666.975	&	6665.858	\\
7	&	1128	&	6835.171	&	6770.489	&	6772.972	&	6775.015	&	6773.924	\\
8	&	1192	&	7221.051	&	7154.631	&	7157.254	&	7159.156	&	7158.158	\\
9	&	1247	&	7552.626	&	7484.752	&	7487.496	&	7489.278	&	7488.356	\\
10	&	1251	&	7576.739	&	7508.761	&	7511.514	&	7513.286	&	7512.371	\\
11	&	1269	&	7685.246	&	7616.801	&	7619.593	&	7621.326	&	7620.435	\\
12	&	1298	&	7860.056	&	7790.865	&	7793.721	&	7795.39	&	7794.539	\\
13	&	1325	&	8022.801	&	7952.925	&	7955.84	&	7957.45	&	7956.635	\\
14	&	1329	&	8046.911	&	7976.934	&	7979.858	&	7981.459	&	7980.649	\\
15	&	1342	&	8125.267	&	8054.962	&	8057.915	&	8059.488	&	8058.695	\\
16	&	1374	&	8318.135	&	8247.033	&	8250.056	&	8251.558	&	8250.809	\\
17	&	1377	&	8336.216	&	8265.04	&	8268.069	&	8269.565	&	8268.819	\\
18	&	1381	&	8360.324	&	8289.049	&	8292.087	&	8293.574	&	8292.833	\\
19	&	1434	&	8679.736	&	8607.166	&	8610.321	&	8611.691	&	8611.019	\\
20	&	1464	&	8860.525	&	8787.232	&	8790.453	&	8791.757	&	8791.124	\\
21	&	1470	&	8896.681	&	8823.245	&	8826.479	&	8827.771	&	8827.145	\\
22	&	1484	&	8981.045	&	8907.276	&	8910.541	&	8911.802	&	8911.194	\\
23	&	1498	&	9065.408	&	8991.307	&	8994.603	&	8995.833	&	8995.242	\\
24	&	1525	&	9228.102	&	9153.367	&	9156.722	&	9157.892	&	9157.336	\\
25	&	1531	&	9264.256	&	9189.38	&	9192.748	&	9193.906	&	9193.356	\\
26	&	1536	&	9294.384	&	9219.391	&	9222.77	&	9223.917	&	9223.373	\\
27	&	1553	&	9396.816	&	9321.429	&	9324.845	&	9325.954	&	9325.432	\\
28	&	1555	&	9408.867	&	9333.433	&	9336.854	&	9337.959	&	9337.439	\\
29	&	1581	&	9565.524	&	9489.491	&	9492.969	&	9494.016	&	9493.528	\\
30	&	1596	&	9655.9	&	9579.524	&	9583.035	&	9584.049	&	9583.58	\\
31	&	1619	&	9794.474	&	9717.575	&	9721.136	&	9722.1	&	9721.658	\\
32	&	1641	&	9927.019	&	9849.624	&	9853.233	&	9854.149	&	9853.733	\\
33	&	1678	&	10149.927	&	10071.706	&	10075.396	&	10076.231	&	10075.859	\\
34	&	1728	&	10451.138	&	10371.816	&	10375.617	&	10376.341	&	10376.028	\\
35	&	1732	&	10475.235	&	10395.825	&	10399.634	&	10400.35	&	10400.042	\\
36	&	1800	&	10884.853	&	10803.976	&	10807.934	&	10808.501	&	10808.27	\\
37	&	1809	&	10939.065	&	10857.996	&	10861.974	&	10862.521	&	10862.3	\\
38	&	1842	&	11137.837	&	11056.069	&	11060.119	&	11060.594	&	11060.41	\\
39	&	1851	&	11192.046	&	11110.089	&	11114.159	&	11114.614	&	11114.44	\\
40	&	1873	&	11324.556	&	11242.137	&	11246.256	&	11246.662	&	11246.514	\\
41	&	1877	&	11348.649	&	11266.146	&	11270.273	&	11270.671	&	11270.527	\\
42	&	1888	&	11414.902	&	11332.17	&	11336.322	&	11336.696	&	11336.563	\\
43	&	1929	&	11661.841	&	11578.261	&	11582.502	&	11582.786	&	11582.699	\\
44	&	1930	&	11667.864	&	11584.263	&	11588.507	&	11588.789	&	11588.703	\\
45	&	1974	&	11932.859	&	11848.361	&	11852.701	&	11852.886	&	11852.848	\\
46	&	1982	&	11981.039	&	11896.379	&	11900.736	&	11900.904	&	11900.874	\\
47	&	1995	&	12059.331	&	11974.407	&	11978.793	&	11978.933	&	11978.917	\\
48	&	2016	&	12185.8	&	12100.454	&	12104.886	&	12104.979	&	12104.986	\\
49	&	2052	&	12402.598	&	12316.534	&	12321.045	&	12321.059	&	12321.104	\\
50	&	2068	&	12498.951	&	12412.569	&	12417.115	&	12417.094	&	12417.156	\\
51	&	2127	&	12854.239	&	12766.7	&	12771.375	&	12771.225	&	12771.348	\\
52	&	2128	&	12860.261	&	12772.702	&	12777.38	&	12777.227	&	12777.351	\\
53	&	2133	&	12890.369	&	12802.713	&	12807.402	&	12807.238	&	12807.368	\\
54	&	2144	&	12956.607	&	12868.737	&	12873.45	&	12873.263	&	12873.403	\\
55	&	2164	&	13077.038	&	12988.782	&	12993.538	&	12993.307	&	12993.468	\\
56	&	2199	&	13287.787	&	13198.859	&	13203.693	&	13203.384	&	13203.581	\\
57	&	2210	&	13354.021	&	13264.884	&	13269.741	&	13269.409	&	13269.617	\\
58	&	2221	&	13420.255	&	13330.908	&	13335.79	&	13335.433	&	13335.652	\\
59	&	2300	&	13895.917	&	13805.083	&	13810.138	&	13809.608	&	13809.906	\\
60	&	2303	&	13913.979	&	13823.09	&	13828.151	&	13827.615	&	13827.916	\\
61	&	2327	&	14058.478	&	13967.143	&	13972.257	&	13971.668	&	13971.993	\\
62	&	2330	&	14076.54	&	13985.15	&	13990.27	&	13989.675	&	13990.002	\\
63	&	2349	&	14190.933	&	14099.192	&	14104.354	&	14103.717	&	14104.063	\\
64	&	2368	&	14305.324	&	14213.234	&	14218.437	&	14217.759	&	14218.124	\\
65	&	2394	&	14461.857	&	14369.291	&	14374.552	&	14373.817	&	14374.207	\\
66	&	2413	&	14576.245	&	14483.334	&	14488.636	&	14487.859	&	14488.267	\\
67	&	2425	&	14648.489	&	14555.36	&	14560.689	&	14559.885	&	14560.305	\\
68	&	2427	&	14660.529	&	14567.365	&	14572.697	&	14571.89	&	14572.312	\\
69	&	2445	&	14768.894	&	14675.405	&	14680.777	&	14679.93	&	14680.369	\\
70	&	2515	&	15190.3	&	15095.56	&	15101.085	&	15100.085	&	15100.59	\\
71	&	2517	&	15202.339	&	15107.564	&	15113.094	&	15112.089	&	15112.597	\\
72	&	2529	&	15274.578	&	15179.591	&	15185.147	&	15184.116	&	15184.635	\\
73	&	2532	&	15292.638	&	15197.597	&	15203.16	&	15202.123	&	15202.644	\\
74	&	2548	&	15388.955	&	15293.633	&	15299.231	&	15298.158	&	15298.694	\\
75	&	2568	&	15509.351	&	15413.677	&	15419.319	&	15418.202	&	15418.757	\\
76	&	2614	&	15786.255	&	15689.779	&	15695.522	&	15694.304	&	15694.902	\\
77	&	2624	&	15846.45	&	15749.801	&	15755.566	&	15754.327	&	15754.933	\\
78	&	2638	&	15930.723	&	15833.833	&	15839.628	&	15838.358	&	15838.977	\\
79	&	2639	&	15936.742	&	15839.835	&	15845.632	&	15844.36	&	15844.98	\\
80	&	2695	&	16273.827	&	16175.959	&	16181.879	&	16180.484	&	16181.156	\\
81	&	2700	&	16303.923	&	16205.97	&	16211.901	&	16210.495	&	16211.171	\\
82	&	2707	&	16346.058	&	16247.986	&	16253.932	&	16252.511	&	16253.193	\\
83	&	2741	&	16550.709	&	16452.061	&	16458.082	&	16456.586	&	16457.299	\\
84	&	2748	&	16592.843	&	16494.077	&	16500.112	&	16498.602	&	16499.321	\\
85	&	2752	&	16616.919	&	16518.085	&	16524.13	&	16522.611	&	16523.333	\\
86	&	2765	&	16695.167	&	16596.114	&	16602.187	&	16600.639	&	16601.374	\\
87	&	2790	&	16845.641	&	16746.17	&	16752.297	&	16750.695	&	16751.452	\\
88	&	2798	&	16893.792	&	16794.187	&	16800.333	&	16798.713	&	16799.476	\\
89	&	2814	&	16990.094	&	16890.223	&	16896.403	&	16894.748	&	16895.526	\\
90	&	2825	&	17056.301	&	16956.247	&	16962.452	&	16960.773	&	16961.56	\\
91	&	2831	&	17092.414	&	16992.261	&	16998.478	&	16996.786	&	16997.579	\\
92	&	2859	&	17260.939	&	17160.323	&	17166.602	&	17164.848	&	17165.665	\\
93	&	2861	&	17272.977	&	17172.327	&	17178.61	&	17176.852	&	17177.672	\\
94	&	2864	&	17291.033	&	17190.334	&	17196.624	&	17194.859	&	17195.681	\\
95	&	2894	&	17471.592	&	17370.401	&	17376.756	&	17374.926	&	17375.774	\\
96	&	2897	&	17489.648	&	17388.407	&	17394.769	&	17392.932	&	17393.783	\\
97	&	2947	&	17790.573	&	17688.518	&	17694.99	&	17693.043	&	17693.937	\\
98	&	2963	&	17886.867	&	17784.554	&	17791.06	&	17789.079	&	17789.986	\\
99	&	2973	&	17947.05	&	17844.576	&	17851.104	&	17849.101	&	17850.017	\\
100	&	2984	&	18013.252	&	17910.6	&	17917.153	&	17915.125	&	17916.051	\\\hline
 \hline
\end{longtable}
\normalsize 

\subsection{Bankruptcy-driven cost allocations}\label{app:numBD}

Table \ref{tab:resultscostbank} presents the corresponding bankruptcy-based cost allocations for the 100 buyers involved in the group-buying situation. For each buyer, the table reports demand, individual cost, separable cost, and the allocations obtained under the four bankruptcy rules considered: $PROP$, $CEA$, $CEL$, and $TAL$.

\footnotesize
\renewcommand{\arraystretch}{0.9}
\begin{longtable}{c|ccc|ccc}
\caption{Demands of product, individual costs, claims and the bankruptcy-driven cost allocations}
\label{tab:resultscostbank}\\
\hline
$i$ & $q_i$ & $C^{GB}(\{i\},q,P)$ & $d_i^{GB}$ &
$\phi_i^{PROP,GB}$ & $\phi_i^{CEA,GB}$ & $\phi_i^{CEL,GB}=\phi_i^{TAL,GB}$\\
\hline
\endfirsthead

\caption[]{Demands of product, individual costs, claims and the bankruptcy-driven cost allocations}\\
\hline
$i$ & $q_i$ & $C^{GB}$ & $d_i^{GB}$ &
$\phi_i^{PROP,GB}$ & $\phi_i^{CEA,GB}$ & $\phi_i^{CEL,GB}=\phi_i^{TAL,GB}$\\
\hline
\endhead
1	&	1021	&	6189.906	&	61.653	&	6131.527	&	6128.253	&	6132.778		\\
2	&	1036	&	6280.374	&	62.088	&	6221.583	&	6218.286	&	6222.811		\\
3	&	1038	&	6292.436	&	62.146	&	6233.590	&	6230.290	&	6234.816		\\
4	&	1039	&	6298.467	&	62.174	&	6239.594	&	6236.293	&	6240.818		\\
5	&	1083	&	6563.818	&	63.428	&	6503.758	&	6500.390	&	6504.915		\\
6	&	1110	&	6726.633	&	64.184	&	6665.858	&	6662.450	&	6666.975		\\
7	&	1128	&	6835.171	&	64.682	&	6773.924	&	6770.489	&	6775.015		\\
8	&	1192	&	7221.051	&	66.420	&	7158.158	&	7154.631	&	7159.156		\\
9	&	1247	&	7552.626	&	67.873	&	7488.356	&	7484.752	&	7489.278	\\
10	&	1251	&	7576.739	&	67.978	&	7512.371	&	7508.761	&	7513.286	\\
11	&	1269	&	7685.246	&	68.445	&	7620.435	&	7616.801	&	7621.326	\\
12	&	1298	&	7860.056	&	69.190	&	7794.539	&	7790.865	&	7795.390	\\
13	&	1325	&	8022.801	&	69.876	&	7956.635	&	7952.925	&	7957.450	\\
14	&	1329	&	8046.911	&	69.977	&	7980.649	&	7976.934	&	7981.459	\\
15	&	1342	&	8125.267	&	70.304	&	8058.695	&	8054.962	&	8059.488 \\
16	&	1374	&	8318.135	&	71.102	&	8250.809	&	8247.033	&	8251.558 \\
17	&	1377	&	8336.216	&	71.176	&	8268.819	&	8265.040	&	8269.565 	\\
18	&	1381	&	8360.324	&	71.275	&	8292.833	&	8289.049	&	8293.574 	\\
19	&	1434	&	8679.736	&	72.571	&	8611.019	&	8607.166	&	8611.691 	\\
20	&	1464	&	8860.525	&	73.292	&	8791.124	&	8787.232	&	8791.757 \\
21	&	1470	&	8896.681	&	73.436	&	8827.145	&	8823.245	&	8827.771 	\\
22	&	1484	&	8981.045	&	73.769	&	8911.194	&	8907.276	&	8911.802 	\\
23	&	1498	&	9065.408	&	74.101	&	8995.242	&	8991.307	&	8995.833 	\\
24	&	1525	&	9228.102	&	74.735	&	9157.336	&	9153.367	&	9157.892	\\
25	&	1531	&	9264.256	&	74.876	&	9193.356	&	9189.380	&	9193.906	\\
26	&	1536	&	9294.384	&	74.992	&	9223.373	&	9219.391	&	9223.917 	\\
27	&	1553	&	9396.816	&	75.387	&	9325.432	&	9321.429	&	9325.954 	\\
28	&	1555	&	9408.867	&	75.433	&	9337.439	&	9333.433	&	9337.959	\\
29	&	1581	&	9565.524	&	76.033	&	9493.528	&	9489.491	&	9494.016		\\
30	&	1596	&	9655.900	&	76.376	&	9583.580	&	9579.524	&	9584.049		\\
31	&	1619	&	9794.474	&	76.899	&	9721.658	&	9717.575	&	9722.100		\\
32	&	1641	&	9927.019	&	77.395	&	9853.733	&	9849.624	&	9854.149		\\
33	&	1678	&	10149.927	&	78.221	&	10075.859	&	10071.706	&	10076.231	\\
34	&	1728	&	10451.138	&	79.322	&	10376.028	&	10371.816	&	10376.341	\\
35	&	1732	&	10475.235	&	79.409	&	10400.042	&	10395.825	&	10400.350		\\
36	&	1800	&	10884.853	&	80.877	&	10808.270	&	10803.976	&	10808.501		\\
37	&	1809	&	10939.065	&	81.069	&	10862.300	&	10857.996	&	10862.521		\\
38	&	1842	&	11137.837	&	81.768	&	11060.410	&	11056.069	&	11060.594		\\
39	&	1851	&	11192.046	&	81.958	&	11114.440	&	11110.089	&	11114.614		\\
40	&	1873	&	11324.556	&	82.419	&	11246.514	&	11242.137	&	11246.662		\\
41	&	1877	&	11348.649	&	82.503	&	11270.527	&	11266.146	&	11270.671		\\
42	&	1888	&	11414.902	&	82.732	&	11336.563	&	11332.170	&	11336.696		\\
43	&	1929	&	11661.841	&	83.580	&	11582.699	&	11578.261	&	11582.786		\\
44	&	1930	&	11667.864	&	83.600	&	11588.703	&	11584.263	&	11588.789		\\
45	&	1974	&	11932.859	&	84.499	&	11852.848	&	11848.361	&	11852.886		\\
46	&	1982	&	11981.039	&	84.661	&	11900.874	&	11896.379	&	11900.904		\\
47	&	1995	&	12059.331	&	84.923	&	11978.917	&	11974.407	&	11978.933		\\
48	&	2016	&	12185.800	&	85.346	&	12104.986	&	12100.454	&	12104.979		\\
49	&	2052	&	12402.598	&	86.064	&	12321.104	&	12316.534	&	12321.059		\\
50	&	2068	&	12498.951	&	86.381	&	12417.156	&	12412.569	&	12417.094		\\
51	&	2127	&	12854.239	&	87.539	&	12771.348	&	12767.699	&	12771.225		\\
52	&	2128	&	12860.261	&	87.559	&	12777.351	&	12773.72	&	12777.227		\\
53	&	2133	&	12890.369	&	87.656	&	12807.368	&	12803.829	&	12807.238		\\
54	&	2144	&	12956.607	&	87.869	&	12873.403	&	12870.066	&	12873.263		\\
55	&	2164	&	13077.038	&	88.256	&	12993.468	&	12990.497	&	12993.307		\\
56	&	2199	&	13287.787	&	88.928	&	13203.581	&	13201.247	&	13203.384		\\
57	&	2210	&	13354.021	&	89.138	&	13269.617	&	13267.481	&	13269.409		\\
58	&	2221	&	13420.255	&	89.347	&	13335.652	&	13333.715	&	13335.433		\\
59	&	2300	&	13895.917	&	90.834	&	13809.906	&	13809.376	&	13809.608	\\
60	&	2303	&	13913.979	&	90.889	&	13827.916	&	13827.439	&	13827.615		\\
61	&	2327	&	14058.478	&	91.335	&	13971.993	&	13971.938	&	13971.668		\\
62	&	2330	&	14076.54	&	91.391	&	13990.002	&	13990.000	&	13989.675		\\
63	&	2349	&	14190.933	&	91.741	&	14104.063	&	14104.393	&	14103.717		\\
64	&	2368	&	14305.324	&	92.09	&	14218.124	&	14218.784	&	14217.759		\\
65	&	2394	&	14461.857	&	92.566	&	14374.207	&	14375.317	&	14373.817		\\
66	&	2413	&	14576.245	&	92.911	&	14488.267	&	14489.704	&	14487.859		\\
67	&	2425	&	14648.489	&	93.128	&	14560.305	&	14561.948	&	14559.885		\\
68	&	2427	&	14660.529	&	93.165	&	14572.312	&	14573.989	&	14571.89		\\
69	&	2445	&	14768.894	&	93.489	&	14680.369	&	14682.354	&	14679.93		\\
70	&	2515	&	15190.300	&	94.740	&	15100.59	&	15103.759	&	15100.085		\\
71	&	2517	&	15202.339	&	94.775	&	15112.597	&	15115.799	&	15112.089		\\
72	&	2529	&	15274.578	&	94.988	&	15184.635	&	15188.038	&	15184.116	\\
73	&	2532	&	15292.638	&	95.041	&	15202.644	&	15206.098	&	15202.123		\\
74	&	2548	&	15388.955	&	95.323	&	15298.694	&	15302.415	&	15298.158		\\
75	&	2568	&	15509.351	&	95.674	&	15418.757	&	15422.811	&	15418.202		\\
76	&	2614	&	15786.255	&	96.475	&	15694.902	&	15699.714	&	15694.304		\\
77	&	2624	&	15846.45	&	96.649	&	15754.933	&	15759.910	&	15754.327		\\
78	&	2638	&	15930.723	&	96.890	&	15838.977	&	15844.183	&	15838.358		\\
79	&	2639	&	15936.742	&	96.908	&	15844.98	&	15850.202	&	15844.360		\\
80	&	2695	&	16273.827	&	97.868	&	16181.156	&	16187.287	&	16180.484		\\
81	&	2700	&	16303.923	&	97.953	&	16211.171	&	16217.383	&	16210.495		\\
82	&	2707	&	16346.058	&	98.072	&	16253.193	&	16259.517	&	16252.511		\\
83	&	2741	&	16550.709	&	98.648	&	16457.299	&	16464.169	&	16456.586		\\
84	&	2748	&	16592.843	&	98.766	&	16499.321	&	16506.303	&	16498.602		\\
85	&	2752	&	16616.919	&	98.834	&	16523.333	&	16530.379	&	16522.611		\\
86	&	2765	&	16695.167	&	99.052	&	16601.374	&	16608.626	&	16600.639		\\
87	&	2790	&	16845.641	&	99.471	&	16751.452	&	16759.101	&	16750.695		\\
88	&	2798	&	16893.792	&	99.605	&	16799.476	&	16807.252	&	16798.713		\\
89	&	2814	&	16990.094	&	99.871	&	16895.526	&	16903.554	&	16894.748		\\
90	&	2825	&	17056.301	&	100.054	&	16961.560	&	16969.761	&	16960.773		\\
91	&	2831	&	17092.414	&	100.154	&	16997.579	&	17005.874	&	16996.786		\\
92	&	2859	&	17260.939	&	100.616	&	17165.665	&	17174.399	&	17164.848		\\
93	&	2861	&	17272.977	&	100.649	&	17177.672	&	17186.436	&	17176.852		\\
94	&	2864	&	17291.033	&	100.699	&	17195.681	&	17204.492	&	17194.859		\\
95	&	2894	&	17471.592	&	101.191	&	17375.774	&	17385.052	&	17374.926		\\
96	&	2897	&	17489.648	&	101.240	&	17393.783	&	17403.107	&	17392.932		\\
97	&	2947	&	17790.573	&	102.054	&	17693.937	&	17704.032	&	17693.043		\\
98	&	2963	&	17886.867	&	102.313	&	17789.986	&	17800.327	&	17789.079		\\
99	&	2973	&	17947.050	&	102.475	&	17850.017	&	17860.51	&	17849.101		\\
100	&	2984	&	18013.252	&	102.652	&	17916.051	&	17926.712	&	17915.125		\\\hline
 \hline
\end{longtable}
\normalsize

\subsection{Variable discount factors}\label{app:discount}

Table \ref{tab:b} presents the variable discount factors associated with the bankruptcy-based cost allocations obtained under the rules $PROP$, $CEA$, $CEL$, and $TAL$.
\footnotesize
\begin{table}[h!]
    \centering\resizebox{0.7\textwidth}{!}{
    \begin{tabular}{c|cccc||c|cccc}\hline
      $i$ & $b^{EP}$ & $b^{PROP}$ & $b^{CEA}$  & $b^{CEL}=b^{TAL}$ &$i$   & $b^{EP}$ & $b^{PROP}$ & $b^{CEA}$  & $b^{CEL}=b^{TAL}$  \\\hline\hline
1	&	0.1408	&	0.1730	&	0.0705	&	0.2121	&	51	&	0.2033	&	0.2027	&	0.1236	&	0.2000	\\
2	&	0.1419	&	0.1735	&	0.0710	&	0.2116	&	52	&	0.2033	&	0.2027	&	0.1240	&	0.2000	\\
3	&	0.1420	&	0.1735	&	0.0711	&	0.2115	&	53	&	0.2036	&	0.2028	&	0.1262	&	0.2000	\\
4	&	0.1421	&	0.1735	&	0.0711	&	0.2115	&	54	&	0.2041	&	0.2031	&	0.1310	&	0.2000	\\
5	&	0.1451	&	0.1750	&	0.0726	&	0.2101	&	55	&	0.2050	&	0.2035	&	0.1397	&	0.2001	\\
6	&	0.1468	&	0.1758	&	0.0735	&	0.2093	&	56	&	0.2067	&	0.2043	&	0.1545	&	0.2001	\\
7	&	0.1480	&	0.1764	&	0.0741	&	0.2089	&	57	&	0.2072	&	0.2046	&	0.1591	&	0.2001	\\
8	&	0.1522	&	0.1784	&	0.0762	&	0.2073	&	58	&	0.2077	&	0.2048	&	0.1637	&	0.2002	\\
9	&	0.1556	&	0.1800	&	0.0779	&	0.2061	&	59	&	0.2114	&	0.2066	&	0.1955	&	0.2003	\\
10	&	0.1559	&	0.1801	&	0.0781	&	0.2060	&	60	&	0.2115	&	0.2066	&	0.1967	&	0.2004	\\
11	&	0.1570	&	0.1807	&	0.0786	&	0.2057	&	61	&	0.2126	&	0.2072	&	0.2060	&	0.2004	\\
12	&	0.1588	&	0.1815	&	0.0795	&	0.2051	&	62	&	0.2128	&	0.2072	&	0.2072	&	0.2004	\\
13	&	0.1604	&	0.1823	&	0.0804	&	0.2047	&	63	&	0.2136	&	0.2076	&	0.2144	&	0.2005	\\
14	&	0.1607	&	0.1824	&	0.0805	&	0.2046	&	64	&	0.2145	&	0.2080	&	0.2216	&	0.2005	\\
15	&	0.1615	&	0.1828	&	0.0809	&	0.2044	&	65	&	0.2157	&	0.2086	&	0.2313	&	0.2006	\\
16	&	0.1634	&	0.1837	&	0.0818	&	0.2039	&	66	&	0.2165	&	0.2090	&	0.2383	&	0.2007	\\
17	&	0.1636	&	0.1838	&	0.0819	&	0.2039	&	67	&	0.2171	&	0.2093	&	0.2426	&	0.2007	\\
18	&	0.1638	&	0.1839	&	0.0820	&	0.2038	&	68	&	0.2171	&	0.2093	&	0.2434	&	0.2007	\\
19	&	0.1669	&	0.1854	&	0.0836	&	0.2031	&	69	&	0.2179	&	0.2097	&	0.2498	&	0.2008	\\
20	&	0.1686	&	0.1862	&	0.0845	&	0.2027	&	70	&	0.2210	&	0.2112	&	0.2744	&	0.2011	\\
21	&	0.1690	&	0.1864	&	0.0846	&	0.2027	&	71	&	0.2211	&	0.2112	&	0.2751	&	0.2011	\\
22	&	0.1698	&	0.1867	&	0.0851	&	0.2025	&	72	&	0.2217	&	0.2115	&	0.2791	&	0.2012	\\
23	&	0.1706	&	0.1871	&	0.0855	&	0.2024	&	73	&	0.2218	&	0.2115	&	0.2802	&	0.2012	\\
24	&	0.1721	&	0.1878	&	0.0862	&	0.2021	&	74	&	0.2225	&	0.2119	&	0.2856	&	0.2012	\\
25	&	0.1725	&	0.1880	&	0.0864	&	0.2020	&	75	&	0.2234	&	0.2123	&	0.2923	&	0.2013	\\
26	&	0.1727	&	0.1881	&	0.0865	&	0.2020	&	76	&	0.2254	&	0.2132	&	0.3074	&	0.2015	\\
27	&	0.1737	&	0.1886	&	0.0870	&	0.2018	&	77	&	0.2258	&	0.2134	&	0.3106	&	0.2016	\\
28	&	0.1738	&	0.1886	&	0.0871	&	0.2018	&	78	&	0.2264	&	0.2137	&	0.3151	&	0.2017	\\
29	&	0.1753	&	0.1893	&	0.0878	&	0.2016	&	79	&	0.2264	&	0.2137	&	0.3154	&	0.2017	\\
30	&	0.1761	&	0.1897	&	0.0882	&	0.2015	&	80	&	0.2288	&	0.2149	&	0.3330	&	0.2020	\\
31	&	0.1774	&	0.1903	&	0.0889	&	0.2013	&	81	&	0.2290	&	0.2150	&	0.3345	&	0.2020	\\
32	&	0.1786	&	0.1909	&	0.0895	&	0.2012	&	82	&	0.2293	&	0.2151	&	0.3367	&	0.2020	\\
33	&	0.1806	&	0.1919	&	0.0905	&	0.2009	&	83	&	0.2308	&	0.2158	&	0.3470	&	0.2022	\\
34	&	0.1832	&	0.1931	&	0.0918	&	0.2007	&	84	&	0.2311	&	0.2160	&	0.3491	&	0.2022	\\
35	&	0.1834	&	0.1932	&	0.0919	&	0.2006	&	85	&	0.2312	&	0.2160	&	0.3503	&	0.2023	\\
36	&	0.1870	&	0.1949	&	0.0937	&	0.2004	&	86	&	0.2318	&	0.2163	&	0.3542	&	0.2023	\\
37	&	0.1875	&	0.1952	&	0.0939	&	0.2003	&	87	&	0.2328	&	0.2168	&	0.3616	&	0.2025	\\
38	&	0.1892	&	0.1960	&	0.0948	&	0.2002	&	88	&	0.2332	&	0.2170	&	0.3640	&	0.2025	\\
39	&	0.1896	&	0.1962	&	0.0950	&	0.2002	&	89	&	0.2338	&	0.2173	&	0.3686	&	0.2026	\\
40	&	0.1908	&	0.1967	&	0.0956	&	0.2002	&	90	&	0.2343	&	0.2175	&	0.3718	&	0.2027	\\
41	&	0.1910	&	0.1968	&	0.0957	&	0.2001	&	91	&	0.2345	&	0.2176	&	0.3735	&	0.2027	\\
42	&	0.1915	&	0.1971	&	0.0960	&	0.2001	&	92	&	0.2357	&	0.2182	&	0.3815	&	0.2029	\\
43	&	0.1936	&	0.1981	&	0.0970	&	0.2001	&	93	&	0.2358	&	0.2182	&	0.3821	&	0.2029	\\
44	&	0.1936	&	0.1981	&	0.0970	&	0.2001	&	94	&	0.2359	&	0.2183	&	0.3829	&	0.2029	\\
45	&	0.1958	&	0.1991	&	0.0982	&	0.2000	&	95	&	0.2371	&	0.2189	&	0.3913	&	0.2031	\\
46	&	0.1962	&	0.1993	&	0.0984	&	0.2000	&	96	&	0.2372	&	0.2189	&	0.3922	&	0.2031	\\
47	&	0.1969	&	0.1996	&	0.0987	&	0.2000	&	97	&	0.2393	&	0.2199	&	0.4059	&	0.2034	\\
48	&	0.1979	&	0.2001	&	0.0992	&	0.2000	&	98	&	0.2399	&	0.2202	&	0.4102	&	0.2035	\\
49	&	0.1997	&	0.2010	&	0.1001	&	0.2000	&	99	&	0.2403	&	0.2204	&	0.4128	&	0.2036	\\
50	&	0.2004	&	0.2013	&	0.1005	&	0.2000	&	100	&	0.2408	&	0.2206	&	0.4158	&	0.2037	\\\hline
    \end{tabular}}
    \caption{Variable discount factors $b^{EP}$,  $b^{PROP}$, $b^{CEA}$ and  $b^{CEL}=b^{TAL}$.}
    \label{tab:b}
\end{table}\normalsize

\section{A simulation study}\label{appD}

The numerical example in Section \ref{sec4}  illustrates how the choice of cost-allocation rule affects the behavior of the variable discount factors across buyers. To assess whether these patterns are robust or specific to the particular numerical instance considered, we conduct a simulation study in which we examine the behavior of both the cost allocations and the associated variable discount factors under a range of scenarios.

For each scenario, we perform 1000 independent simulations, considering three group sizes, $|N|\in\{10,100,1000\}$, and three parameter regimes $0<\eta<1/2$, $\eta=1/2$, and $1/2<\eta<1$. The quantities $q$ and the parameters $a$ and $b$ are generated according to the same simulation procedure and assumptions as in the preceding numerical example. For each replication, we classify the resulting sequence across buyers as increasing (I), decreasing (D), or non-monotonic (NM). We report the percentage of replications falling into each category, separately for each group size, parameter regime, and allocation rule.

Tables \ref{tab:simulationstudy:allocations_rules} and \ref{tab:simulationstudy:vdf} summarize the results for the cost allocations and the variable discount factors, respectively.

\footnotesize
\begin{table}[h!]
    \centering\resizebox{0.8\textwidth}{!}{
    \begin{tabular}{c|c|c||cccc}\hline
 &  &  & \multicolumn{4}{c}{Bankruptcy-driven cost allocations $\phi^\chi(N,q,P)$}
         \\\hline
$\eta$ &    $|N|$ & $EP$  &  $PROP$ & $CEA$ & $CEL$&$TAL$     \\\hline\hline
$\eta<1/2$&$10$& (100/0/0)& (100/0/0) & (100/0/0) &  (100/0/0) & (100/0/0)  \\
&$100$& (100/0/0)& (100/0/0) & (100/0/0) &  (100/0/0) & (100/0/0)  \\
&$1000$&(100/0/0)& (100/0/0) & (100/0/0) &  (100/0/0) & (100/0/0)  \\
\hline
$\eta=1/2$&$10$&(100/0/0)& (100/0/0) & (100/0/0) &  (100/0/0) & (100/0/0)  \\
&$100$& (100/0/0)& (100/0/0) & (100/0/0) &  (100/0/0) & (100/0/0)  \\
&$1000$&(100/0/0)& (100/0/0) & (100/0/0) &  (100/0/0) & (100/0/0)  \\\hline
$\eta>1/2$&$10$& (100/0/0)& (100/0/0) & (100/0/0) &  (100/0/0) & (100/0/0)  \\
&$100$& (100/0/0)& (100/0/0) & (100/0/0) &  (100/0/0) & (100/0/0)  \\
&$1000$&(100/0/0)& (100/0/0) & (100/0/0) &  (100/0/0) & (100/0/0)  \\
\hline
    \end{tabular}}
    \caption{Percentage of simulations exhibiting increasing (I),
    decreasing (D), and non-monotonic (NM) behavior of bankruptcy-driven cost allocations for different values of $n$.Each entry reports the percentage of the 1000 replications in the order $(I/D/NM)$, corresponding to increasing, decreasing, and non-monotonic behavior, respectively.}
    \label{tab:simulationstudy:allocations_rules}
\end{table}\normalsize

\footnotesize
\begin{table}[h!]
    \centering\resizebox{0.8\textwidth}{!}{
    \begin{tabular}{c|c|c||cccc}\hline
    
 &   &    & \multicolumn{4}{c}{Variable discount factor $b^\chi(q_i,q_{-i})$} \\\hline
  &  $|N|$ & $EP$ &  $PROP$ & $CEA$ & $CEL$&$TAL$     \\\hline\hline
    $\eta<1/2$&$10$& (100/0/0)& (100/0/0) & (100/0/0) &  (0/0/100) & (0/0/100)  \\
&$100$& (100/0/0)& (100/0/0) & (100/0/0) &  (0/0/100) & (0/0/100)  \\
&$1000$& (100/0/0)& (100/0/0) & (100/0/0) &  (0/0/100) & (0/0/100)  \\
\hline
$\eta=1/2$& $10$& (100/0/0)& (100/0/0) & (100/0/0) &  (0/0/100) & (0/0/100)  \\
&$100$ &  (100/0/0)& (100/0/0) & (100/0/0) &  (0/0/100) & (0/0/100)  \\
&$1000$ & (100/0/0)& (100/0/0) & (100/0/0) &  (0/0/100) & (0/0/100)  \\
\hline$\eta>1/2$&$10$& (100/0/0)& (100/0/0) & (100/0/0) &  (0/0/100) & (0/0/100)  \\
&$100$& (100/0/0)& (100/0/0) & (100/0/0) &  (0/0/100) & (0/0/100)  \\
&$1000$& (100/0/0)& (100/0/0) & (100/0/0) &  (0/0/100) & (0/0/100)  \\
\hline
    \end{tabular}}
    \caption{Percentage of simulations exhibiting increasing (I),
    decreasing (D), and non-monotonic (NM) behavior of 
    variable discount factor  for different values of $n$.Each entry reports the percentage of the 1000 replications in the order $(I/D/NM)$, corresponding to increasing, decreasing, and non-monotonic behavior, respectively.}
    \label{tab:simulationstudy:vdf}
\end{table}\normalsize

The simulation results are fully consistent with the patterns observed in the numerical example in Section \ref{sec4}. In particular, the cost allocations exhibit increasing behavior across buyers under all allocation rules, regardless of the group size or the parameter regime considered. The variable discount factors display a more differentiated pattern: they are increasing under the EP, PROP, and CEA rules, whereas they exhibit non-monotonic behavior under the CEL and TAL rules. These patterns remain unchanged across all group sizes and parameter regimes considered. Overall, the simulation study confirms that the qualitative findings illustrated by the numerical example are robust to changes in group size and parameter values.

\end{document}